\documentclass[12pt]{article}
\usepackage{amsmath,amssymb,amscd,latexsym,mathrsfs}
\usepackage[all]{xy}
\usepackage[cp1251]{inputenc} %
\usepackage[english]{babel} %

\newcommand{\coLim}{\underrightarrow{\lim}}
\newcommand{\Lim}{\underleftarrow{\lim}}
\newcommand{\dom}{{\rm dom\,}}
\newcommand{\cod}{{\rm cod\,}}
\newcommand{\Imm}{{\rm Im\,}}
\newcommand{\Ker}{{\rm Ker\,}}
\newcommand{\Coker}{{\rm Coker\,}}

\newcommand{\Set}{{\rm Set}}
\newcommand{\Cat}{{\rm Cat}}
\newcommand{\Ab}{{\rm Ab}}
\newcommand{\Hot}{{\rm Hot}}
\newcommand{\Ob}{{\rm Ob}}

\newcommand{\pt}{{\rm pt}}
\newcommand{\nr}{{\rm B}}
\newcommand{\nd}{{\rm N}}

\newcommand{\NN}{{\,\mathbb N}}

\newcommand{\mC}{{\,\mathscr C}}
\newcommand{\mD}{{\,\mathscr D}}

\newcommand{\RR}{{\,\mathbb R}}

\newcommand{\ZZ}{{\,\mathbb Z}}
\newcommand{\II}{{\,\mathbb I}}
\newcommand{\mA}{{\mathcal A}}

\newcommand{\fF}{{\frak F}}

\newtheorem{theorem}{\bf Theorem}[section]
\newtheorem{lemma}[theorem]{\bf Lemma}
\newtheorem{proposition}[theorem]{\bf Proposition}
\newtheorem{corollary}[theorem]{\bf Corollary}
\newtheorem{definition}{\sc Definition}[section]
\newtheorem{example}[definition]{\sc Example}
\newtheorem{remark}[definition]{\sc Remark}

\def\leq{\leqslant}
\def\geq{\geqslant}

\begin{document}

\begin{center}
{\large Homology and cohomology of cubical sets\\
 with coefficients in systems of objects}\\ 
 \medskip
{Ahmet A. Husainov}
\end{center}

\begin{abstract}
This paper continues the research of the author on the homology of cubical and semi-cubical sets with coefficients in systems.  The main result is the theorem that the homology of cubical sets with coefficients in contravariant systems in an Abelian category with exact coproducts is isomorphic to the left satellites of a colimit functor. 
It allowed proving a number of new assertions presented in this paper about homology and cohomology of cubical sets with coefficients in systems, including homology and cohomology with coefficients in local systems.
\end{abstract}

2010 Mathematics Subject Classification 55N25, 55U05, 18G35, 18A40, 18E10, 18G10, 18G40, 55P10, 55U15


\tableofcontents

\section{Introduction}
The work is devoted to the homology of cubical sets with coefficients in contravariant systems 
of objects in Abelian category with exact coproduct (shortly AB4-category).

This paper was preceded by the author's papers on the homology of semi-cubical sets \cite{X2008} and cubical sets \cite{X2019} with coefficients in contravariant systems of abelian groups.

Outlining the history and applications of cubical homology, we consider the homology theory of semi-cubical sets with coefficients in systems as part of the homology theory  of cubical sets, although this is proved later (Theorem \ref{semicubecube}).

We advise the reader to start acquaintance with the article from the preliminary section, and to read the introduction only after reading the text of the entire article.

\subsection{Brief history and applications of cubical homology}

Serre \cite{ser1951} introduced cubical singular homology of topological spaces with coefficients in local systems and applied them to construct spectral sequences used in the calculation of homotopy groups.
Eilenberg and MacLane \cite{eil1953} proved the equivalence of the cubical singular
and simplicial singular homology 
theories. Abstract cubical sets and their homology and homotopy groups were introduced by Kan \cite{kan1955}.

The theory of cubical sets led to the development of non-Abelian algebraic topology \cite{bro2011}. One of the most important structures, the $\omega$-groupoid, was introduced in \cite{bro1981} as a cubical set with additional operations. 

To construct invariants of groups and knots, rack cohomology and quandle cohomology were used, allowing cubical description \cite{cov2012}, \cite{cov2014}, \cite{cov2019}.
Classical methods from \cite{eil1953} were also used for non-standard homologies. For example, to study solutions of the Yang-Baxter equation, skew cubical homology in the category of $R$-modules \cite{leb2017} is considered.

A cubical set is defined as a presheaf of sets over some small category called a cubical site.
Grandis and Mauri \cite{gra2003} drew attention to the most popular cubical sites, proved the existence of a normal form of morphism decomposition for each of these sites, and characterized an extended cubical site (with connections and interchanges).
Cubical homology with connections was also studied in \cite{bar2021}.

  Cubical discrete homology \cite{bar2014}, \cite{bar2019}, \cite{bar2020} appeared to study metric spaces.

The homology of semi-cubical sets was used to study the topological properties of models for parallel computing systems - automata of higher dimensions \cite{gau2003}, \cite{fah2005},
\cite{faj2006}, \cite{faj2016}, \cite{kah2018}.
To model parallel computing processes, Grandis developed a directed algebraic topology \cite{gra2009}, in which the directed cubical homology \cite{gra2005} introduced by him found applications. 

The homology of semi-cubical sets allowed the author of this article to solve the problem posed by him in 2004 and consisting in constructing an algorithm for calculating the homology groups of safe Petri nets \cite{X2013}.

The singular cohomology theory is used to formalize the synthetic cohomology theory in the cubical extension of Agda \cite{bru2022}. Synthetic cohomology groups are used to check that topological spaces are not homotopy equivalent. There is a problem related to the choice of a basis for singular cochains. The operation of choosing this basis is not invariant under
 homotopy equivalence, which makes this operation impossible in the synthetic
  formalization
cohomology groups. But the developers found a way out using the Eilenberg-McLane spaces. These cohomologies were used to formalize cellular cohomologies \cite{buc2018}.
The formalization of the synthetic cohomology theory in Agda made it possible to simplify many proofs from the theory of homotopy types \cite{bru2022}. The optimized definition of the cup product made it possible to give a complete system of axioms for transforming integer cohomology groups into a graded commutative ring. This has been used to characterize the cohomology groups of spheres, the torus, the Klein bottle, and the real/projective planes. Constructive proofs allow Cubical Agda to be used to distinguish between spaces through computation.

Cubical homology is used in digital topology, for image analysis, and for analyzing topological data \cite{kac2004}. Pixels, voxels and similar dot patterns are represented as elementary cubes. They can form a semi-cubical set \cite{X2008}. One of the most important tasks is the construction of an algorithm for detecting simple points - points, the removal of which does not change the topology of the processed figure.
This problem is solved in \cite{nie2006} using cubical homology. Cubical
 homology is also applied to the algorithmic calculation of the Conley index for continuous mappings \cite{pil2008}.
The paper \cite{jam2020} develops numerical cubical homology and proves Gurevich's theorem for numerical cubical singular homology.
The work \cite{jam2022} develops the theory of numerical cubical homology for the analysis of digital images.

To process large volumes of topological data, it becomes necessary to calculate invariants with respect to the similarity transformation.
For example, for fractal sets, it may be necessary to calculate fractional Betti numbers. This is done by passing to the limit, similar to that used to calculate the Hausdorff box dimension. This problem is solved using persistent homology \cite{car2009}, \cite{ede2002}, and simplicial homology is traditionally used.
To solve similar problems, when calculating persistent homology in a metric space, \cite{cho2021} used cubical homology. Calculation algorithms are more efficient than algorithms
based on triangulation methods, and they can naturally be used for sets of points in an $n$-dimensional space with the distance between points calculated using the $L_{\infty}$ norm.
Similar questions were studied in \cite{gra20032} for image analysis, both with the help of simplicial and cubical homology.

\subsection{The outline of the paper}

\subsubsection{Definition of homology for a cubical set}

For a cubical object in the Abelian category $\mA$ 
given as a functor $F: \Box^{op}\to \mA$ the normalized complex $C^N_*(F)$
is defined. It consists of objects such that,
for each $k\geq 0$, the object $C^N_k(F)$ is the quotient object of $F(\II^k)$ modulo 
the subobject 
of degenerate chains. So
$$
C^{\nd}_k(F)= \Coker (F(\II^{k-1})^{\oplus k} 
\xrightarrow{(F(\sigma^k_1), \ldots, F(\sigma^k_k))}
 F(\II^k)),
$$
where $\sigma^k_i: \II^k\to \II^{k-1}$ are the degeneracy morphisms  
in the category  $\Box$ of cubes.
The differentials $d_k$ of $C^N_*(F)$ are induced by morphisms
$$
\sum^k_{i=1}(-1)^i(F(\delta^{k,0}_i) -F(\delta^{k,1}_i)): 
F(\II^k)\to F(\II^{k-1}),
$$
for all $k\geq 1$, and the differential $d_0=0$ is added.
This normalized complex for a cubical object 
in Abelian category was explored in  \`Swi\c{a}tek 
\cite{swi1981}.

Homology objects 
$H^N_k(F):= \Ker d_k / \Imm d_{k+1}$ of the compex $C^N_*(F)$ are called 
the $k$-th homology objects (or simply homology) of a cubical object $F$, 
for all $k\geq 0$.

The homology of cubical set $X: \Box^{op}\to \Set$
with coefficients in the functor  $G: (\Box/X)^{op}\to \mA$
is defined as the homology of a cubical object equal to the left Kan extension
$Lan^{Q^{op}_X}G$ of the functor $G: (\Box/X)^{op}\to \mA$ along the functor
 $Q^{op}_X: (\Box/X)^{op} \to \Box^{op}$,
$$
H_k(X, G) = H^N_k( Lan^{Q^{op}_X}G), k\geq 0.
$$

\subsubsection{Main results}

In \cite{X2019} it is proved that the homology groups $H_n(X,G)$ of the cubical set $X$ with coefficients in the contravariant system of Abelian groups $G: (\Box/X)^{op} \to \Ab $ are isomorphic to values on $G$ of the left $n$th derived of the colimit functor 
$\coLim^{(\Box/X)^{op}}_n G$ for all $n\geq 0$.

This statement made it possible to apply the homology theory of small categories 
to the study of homology of cubical sets.
Unfortunately, the principle of duality (in category theory) applied to this statement does not give 
a similar result for cohomology (it leads to an isomorphism of cohomology of cubical sets and 
derived of the limit functor for covariant systems of compact abelian
groups).

We prove that the homology of the cubical set $X$ with coefficients in the contravariant systems $G: (\Box/X)^{op} \to \mA$ in the AB4-category $\mA$ are isomorphic to the left satellites of the colimit functor
$\coLim^{(\Box/X)^{op}}_n G$ (Theorem \ref{main2}).
The principle of duality gives an isomorphism between the
cohomology of a cubical set $H^n(X, G)$ and the values of
 right satellites of the limit functor $\Lim^n_{\Box/X}G$
 with coefficients in the covariant system $G: \Box/X\to \mA$
  in the AB4*-category $\mA$ (Corollary \ref{main2coh}).
  
  These results are preceded by an auxiliary Theorem \ref{main1} on the homology of a cubical object.
The proofs of all other assertions use Theorem \ref{main2}.

A theorem on the invariance of homology of a cubical set with coefficients in a contravariant system with respect to the transition to a direct image is proved.

A criterion for the invariance of the homology of a cubical set with respect to the transition to the inverse image for contravariant systems with values in arbitrary AB4-categories has been supplemented.

It is proved that the Baues-Wirshing cohomology of a small category with coefficients in the natural system is isomorphic to the cohomology of the cubical nerve of this category with coefficients in the covariant system corresponding to the natural system.

Complexes and formulas for homology of cubical sets with coefficients in local systems are obtained.
The isomorphism of the cohomology in local systems on weak equivalent cubical sets with respect to the standard model structure is proved.

A spectral sequence is constructed that converges to the cohomology of the colimit of 
cubical sets with covariant systems.

For the cohomology of cubical sets with coefficients in local systems of abelian groups, a spectral sequence of a morphism of cubical sets is constructed whose inverse fibers are weak
 equivalent to each other (with respect to the standard model structure).

It is proved that the homology of the semi-cubical set with coefficients in the contravariant system is isomorphic to the homology of the universal cubical set with coefficients in the extended contravariant system.

\subsubsection{Presentation steps}

Section 1 contains an introduction, a brief history and applications of cubical homology, and a description of the results obtained.
 
Section 2 is devoted to preliminary information.
Notation is introduced and a definition of homology of small categories with coefficients in a diagram of objects in an AB4-category is given.
A method is given for constructing a complex whose homology is equal to that of the category. The method is justified in Proposition 2.3.

A $\mD$-set is a functor $\mD^{op}\to \Set$, where $\mD$ is an arbitrary small category.

Section 3 is devoted to the homology of $\mD$-sets $X$ with coefficients in the contravariant system $(\mD/X)^{op}\to \mA$.
The homology of $\mD$-sets is the sequence of left satellites of the colimit functor
$\coLim^{(\mD/X)^{op}}: \mA^{(\mD/X)^{op}}\to \mA$, where $\mD$ is an arbitrary small category, $\mA$ is an AB4-category , $X$ is a diagram of sets over 
$\mD^{op}$.

Formulas are given that reduce the calculation of the homology of an $\mD$-set with coefficients in contravariant systems to the calculation of the homology of the category $\mD^{op}$ (proposition \ref{contralan}).
 A corollary containing a dual assertion to the proposition \ref{contralan} is given for covariant systems.
 Examples are given showing how these formulas are arranged for the case of contravariant and covariant systems of abelian groups.
A method for proving the assertion \ref{contralan} using discrete Grothendieck fibrations is mentioned (Remark \ref{fibrgroth}).
 
  For an arbitrary morphism of set diagrams $f: X\to Y$ and a contravariant system $F$ on $X$, formulas are given for transforming this contravariant system when passing to the direct image of $f_* F$ on $Y$ (proposition \ref{landiscr} ) and proved
 that the homology of $Y$ with coefficients in $f_* F$ is isomorphic to the homology of $X$ with coefficients in $F$ (Theorem \ref{lanhomol}).
 
 A criterion for the invariance of the homology of a $\mD$-set with coefficients in a contravariant system in $\mA$ under the inversion of e for each Abelian category $\mA$ is given (proposition \ref{critiso}).
 We present a spectral sequence for the colimit homology of diagrams of $\mD$-sets, constructed in \cite{X1989}, and a generalized covering spectral sequence for a morphism of $\mD$-sets, constructed in \cite{X1991}.
 
 Section 4 constructs the projective resolution of the cocubical object $\Delta_{\Box}\ZZ$ 
 in the category $\Ab^{\Box}$. For any cubical object $F\in \mA^{\Box^{op}}$ 
 of the AB4-category $\mA$, the tensor product of this resolution
   and the cubical object $F$ gives a complex whose homology is equal to 
   $\coLim^{\Box^{op}}_n F$, for all $n\geq 0$.
It is proved that this complex is isomorphic to the normalized complex of the cubical
 object $F$, which implies that the normalized $n$th homology of the cubical object $F$ is isomorphic to $\coLim^{\Box^{op}}_n F$ (Theorem \ref{main1}).
An example is given showing what a normalized complex looks like for
 a cubical Abelian group.
 
 Section 5 is devoted to the theorem on the isomorphism of (normalized) homologies of 
 a cubical set with coefficients in contravariant systems to the homology of the category of its cubes, as well as the consequences of this theorem.
First, formulas are constructed for obtaining a normalized complex for a contravariant 
system of objects on a cubical set.
An example of the construction of this complex for a contravariant system of Abelian groups is given. It is proved that the $n$th homology of the normalized complex for the contravariant system $G$ on the cubical set $X$ is isomorphic to $\coLim^{(\Box/X)^{op}}_n G$ (Theorem \ref{main2} ).
A similar assertion has been proved for cohomologies of cubical sets with coefficients 
in covariant systems (Corollary \ref{main2coh}).

The invariance of the homology of the cubical set with coefficients in the contravariant system is obtained under the transition to the forward image with respect to the morphism 
of cubical sets (Corollary \ref{dirhomol}) and the criterion for homology invariance under the transition to the inverse image with respect to the morphism of cubical sets (Corollary \ref{critcube}).

The spectral sequence of the colimit of $\mD$-sets leads to the spectral sequence
 of the colimit of cubical sets with contravariant systems (Corolary \ref{coveringcor}).

At the beginning of the article, we introduced the homology of small categories from the book by Gabriel and Zisman \cite{gab1967} as simplicial homology.
Corollary \ref{homolcatcub} shows that the homology of small categories is 
isomorphic to cubical homology.
Baues and Wirsching \cite{bau1985} introduced the cohomology of small categories with coefficients in natural systems as cohomology of cosimplicial abelian groups. Corollary \ref{homolbwcub} shows that this cohomology can be regarded as cubical cohomology.
 
Section 6 investigates the homology and cohomology of cubical sets with coefficients in local systems.
A complex for homology with coefficients in local systems is obtained, in which chains are equal to zero under degenerate cubes (Theorem \ref{comloc}).

A short introduction to the Grothendieck test categories \cite{gro1983} is given and the definition of a standard model structure is recalled.
It is proved that for a weak equivalence with respect to the standard model structure between the cubical sets $f: X\to Y$ and the local system of objects of Abelian groups $L: \Box/Y \to \Ab$ there exists a natural isomorphism $H^n(Y, L )\to H^n(X, f^*L)$ cohomology groups $Y$ with coefficients in $L$ and homology groups $X$ with coefficients in the local system 
$f^* L := L\circ \Box/f $  (Corollary \ref{lociso}).

For a cubical set morphism $f: X\to Y$ whose diagram of inverse fibers consists of weak equivalences and a local system of Abelian groups $G$ on $X$, 
we construct a spectral sequence connecting
on the one hand, the cohomology of $Y$ with coefficients in the local system of cohomology of inverse fibers of the morphism $f$ over $Y$ and, on the other hand, the cohomology of the cubical set $X$ with coefficients in $G$ (Corollary \ref{spseqserr}).

 Section 7 is devoted to the homology of a semi-cubical set with coefficients in a system of objects in the AB4-category.
A universal cubical set is constructed, into which a semi-cubical set is embedded. It is shown that any contravariant system on a semi-cubical set can be extended to this universal cubical set, and the homology with coefficients in the contravariant system are isomorphic to the original one.

Let $X$ be a semi-cubical set. The universal cubical set containing it is constructed as a left Kan extension of the functor $X$ along the embedding $J^{op}: \Box^{op}_+ \subset \Box^{op}$.
This cubical set consists of the coproducts $(Lan^{J^{op}}X)_n= \coprod\nolimits_{\II^n\stackrel{\gamma} \twoheadrightarrow\II^k}X_k$ (Proposition \ref{cubeforsemi}). It is proved that the embedding functor of the category $\Box_+/X$ into the category $\Box/ Lan^{J^{op}}X$ has a left adjoint (Proposition \ref{semcubadj}).
Whence it follows that the homology of the universal set
with coefficients in the composition $(\Box/Lan^{J^{op}} X)^{op} \xrightarrow{S^{op}} (\Box_+/X)^{op} \xrightarrow{ F} \mA$ are isomorphic to the homology of the semi-cubical set $X$ with coefficients in $F$ (Theorem \ref{semicubecube}).

\section{Preliminaries}

\subsection{Notation}

The following notation will apply:

$\Set$  - category of sets and mappings;

$\Ab$  - category of abelian groups and homomorphisms;

$\pt$ - a category consisting of a single object
and a unique morphism;

$in$ - morphism of the coproduct cone or monomorphism into an object;

$pr$ - morphism of the product cone or epimorphism onto an object;

$\ZZ(-): \Set \rightarrow \Ab$ is a functor assigning to each set $E$ a free abelian group $\ZZ(E)$ with basis $E$ and to each mapping $f: E_1\rightarrow E_2$ - the canonical homomorphism $\ZZ(f): \ZZ(E_1)\rightarrow \ZZ(E_2)$ extending this mapping;

$\II$ - linearly ordered set $\{0,1\}$ with the smallest order relation, containing the pair $0\leq 1$;

$\ZZ$ - set or additive group of integers;

$\NN$ - set of non-negative integers;

$\RR$ - set of real numbers;

$\Delta$ - category whose objects are final
linearly ordered sets $[n]=\{0, 1, \ldots, n\}$, where $n\in \NN$, and sets of morphisms $[m]\to [n]$, for all $m, n\in \NN$, consist of non-decreasing mappings.

For an arbitrary category $\mA$, let $\mA^{op}$ denote the dual category. For objects $a,b\in \mA$, let $\mA(a,b)$ denote the set of morphisms $a\rightarrow b$.

If $\mC$ is a small category, then a functor $\mC\to \mA$ is called a \textit{diagram of objects} 
in the category $\mA$.
In some cases, diagrams will be conveniently denoted by
 $\{X^c\}_{c\in \mC}$, or $\{X^c\}$ for short, indicating the values $X^c$ of this functor on the objects $c\in \Ob(\mC )$.
The categories of diagrams and natural transformations between them will be denoted by $\mA^{\mC}$.
A presheaf of objects of the category $\mA$ in $\mC$ or a $\mC$-object on $\mA$ is a contravariant functor from $\mC$ to $\mA$.

Let $S: \mC\to \mD$ be a functor between small categories.
For an arbitrary category $\mA$, by $(-)\circ S: \mA^{\mD}\to \mA^{\mC}$ we denote the functor carrying each functor $F: \mD\to \mA $ into the functor $F\circ S: \mC\to \mA$, and the natural transformation $\eta: F \to G $ into $\eta*S: FS\to GS$. Here $\eta*S$ is the natural transformation
defined as $(\eta*S)_c= \eta_{S(c)}: FS(c)\to GS(c)$ for all $c\in \Ob\mC$.
The left adjoint functor to $(-)\circ S$, the functor of the left Kan extension \cite{mac1972}, is denoted by $Lan^S: \mA^{\mC}\to \mA^{\mD}$, and the right adjoint to $(-)\circ S$, the right Kan extension  \cite{mac1972}, is denoted by $Ran_S$.
For any object $d\in \mD$ we call the left fiber $S$ over $d$ (respectively, the right fiber $S$ under $d$) and denote by $S/d$ (respectively $d/S$) the category objects over $d$ (respectively, under $d$) in the sense of \cite[Page 45]{mac1972}. 
If $S: \mC\xrightarrow{\subseteq} \mD$ is a full embedding of a subcategory, then the category $S/d$ (respectively $d/S$) is denoted by $\mC/d$ (respectively $d/\mC$ ).

$\Delta_{\mC}\ZZ$ or $\Delta\ZZ$ denotes the diagram of Abelian groups $\mC\to \Ab$ taking constant values $\ZZ$ on objects, and $1_{\ZZ}$ - on morphisms.

By an AB4-category (or an Abelian category with exact coproducts) we mean an Abelian category satisfying axioms AB3-AB4 from Grothendieck's book \cite{gro1957}. Dually, an AB4*-category (or an Abelian category with exact products) is an Abelian category satisfying the axioms AB3*-AB4*.
For an AB4-category $\mA$ and the small category $\mC$ we consider the left satellites of the right exact additive colimit functor $\coLim^{\mC}: \mA^{\mC}\to \mA$ defined in all non-negative dimensions, in the sense of 
\cite[\S2.2]{gro1957}. The values of this satellite on objects of the category $\mA^{\mC}$ are described in \cite[Appendix II]{gab1967}.

By simplicial sets we mean the functors $\Delta^{op}\to \Set$. Morphisms of simplicial sets are natural transformations.

Objects of the category $\Delta$ can be considered as categories, and morphisms $[m]\to [n]$ as functors.

Let $\mC$ be a small category. Its {\it nerve} $\nr\mC: \Delta^{op}\to \Set$ is a simplicial set that assigns to each $[n]\in \Delta$ a set $\nr_n\mC$ of functors $[n] \to\mC$, and the morphism $[m]\to [n]$ is the mapping $\nr_n\mC\to \nr_m\mC$ acting on elements $[n]\to\mC$ from $\nr_n \mC$ as a composition $[m]\to[n]\to \mC$.

\subsection{Homology of small categories}

We consider the homology of small categories with coefficients in diagrams
 of objects in AB4-categories, as the homology of the complex
 of objects constructed in the book \cite{gab1967}. Then we indicate a method for constructing other complexes whose homology objects are naturally isomorphic to the homology of small categories with coefficients in these diagrams. We will need this method to study homology
 of the category of cubes

\subsubsection{A chain compex for the homology of small categories}

Let us introduce the homology of a small category with coefficients in a 
diagram of objects.

\begin{definition}
Let $\mC$ be a small category and let $\mA$ be an AB4-category.
 The homology functors of the category $\mC$ are the left satellites $\coLim^{\mC}_n: \mA^{\mC}\to \mA$ (in the sense of \cite{gro1957}) of the colimit functor $\coLim^{\mC }: \mA^{\mC}\to \mA$. Its values on $F\in \mA^{\mC}$ are called the ($n$-th) homology of the category $\mC$ with coefficients in $F$.
 \end{definition}

Construct a complex for calculating the homology of the category $\mC$ with coefficients in $F: \mC\to \mA$.
For each $n\geq 0$ we denote by
$$
C_n(\mC,F)= \bigoplus\limits_{c_0\stackrel{\alpha_1}\to
\cdots \stackrel{\alpha_n}\to c_n}F(c_0)
$$
an object equal to the coproduct of a family of objects whose indices run through all sequences of composable morphisms of length $n$.
Each sequence of $n$ composable morphisms
$s= (c_0\stackrel{\alpha_1}\to \cdots \stackrel{\alpha_n}\to c_n)$ matches
$(n-1)$-sequences
$$
\delta^n_i(s) =
\begin{cases}
c_1\stackrel{\alpha_2}\to \cdots \stackrel{\alpha_n}\to c_n, & \text{if $i=0$;} \\
c_0\stackrel{\alpha_1}\to \cdots \to c_{i-1}\stackrel{\alpha_{i+1}\alpha_i}\to c_{i+1}
\to \cdots \stackrel{\alpha_n}\to c_n, & \text{if $0<i<n$;}\\
c_0\stackrel{\alpha_1}\to \cdots \stackrel{\alpha_{n-1}}\to c_{n-1}, & \text{if $i=n$}.
\end{cases}
$$
Denote by $\lambda_{c_0\stackrel{\alpha_1}\to
\cdots \to \stackrel{\alpha_n}\to c_n}: F(c_0)\to
\bigoplus\limits_{c_0\stackrel{\alpha_1}\to
\cdots \to \stackrel{\alpha_n}\to c_n}F(c_0)$ 
the morphisms of the coproduct cone.

For each $i$ in the range $1\leq i\leq n$ there is a unique morphism $d_n^i: C_n(\mC,F)\to C_{n-1}(\mC,F)$ making the following diagram commutative
$$
\xymatrix{ C_n(\mC,F) \ar[rr]^{d_n^i} && C_{n-1}(\mC,F)\\
& F(c_0) \ar[lu]^{\lambda_{c_0\stackrel{\alpha_1}\to
\cdots \to \stackrel{\alpha_n}\to c_n}} \ar[ru]_{\quad\lambda_{\delta^n_i(c_0\stackrel{\alpha_1}\to
\cdots \to \stackrel{\alpha_n}\to c_n)}}
}
$$
There is also a unique morphism $d_n^0: C_n(\mC,F)\to C_{n-1}(\mC,F)$ which makes the diagram commutative:
$$
\xymatrix{ C_n(\mC,F) \ar[rr]^{d_n^i} && C_{n-1}(\mC,F)\\
F(c_0) \ar[u]^(.45){\lambda_{c_0\stackrel{\alpha_1}\to
\cdots \stackrel{\alpha_n}\to c_n}} \ar[rr]_{F(\alpha_1)} &&
F(c_1) \ar[u]_(.45){\lambda_{\delta^n_0(c_0\stackrel{\alpha_1}\to
\cdots \stackrel{\alpha_n}\to c_n)}}
}
$$

As a result, we obtain a chain complex of objects and morphisms of the category $\mA$:
$$
 0\leftarrow C_0(\mC,F) \stackrel{d_1}\leftarrow C_1(\mC,F) \stackrel{d_2}\leftarrow
 \cdots \leftarrow C_{n-1}(\mC,F) \stackrel{d_n}\leftarrow C_n(\mC,F) \leftarrow \cdots,
$$
whose differentials are defined by the formula
$d_n= \sum\limits_{i=0}^n (-1)^i d_n^i$. Denote this complex by
$C_*(\mC,F)$.
The correspondence $F\mapsto C_*(\mC,F)$ gives a functor from the category  $\mA^{\mC}$
to the category of chain complexes in the category $\mA$.

\begin{proposition}\cite[Appendix II, Proposition 3.3]{gab1967}
If an Abelian category $\mA$ has exact coproducts, then for every small category $\mC$ there exists a unique sequence of left satellites $\coLim^{\mC}_n: \mA^{\mC}\to \mA$, $n \geq 0$, of the colimit functor
$\coLim^{\mC}: \mA^{\mC}\to \mA$.
The values of these satellites on $F\in \mA^{\mC}$ are naturally isomorphic to the objects
homology of the chain complex $C_*(\mC,F)$.\end{proposition}

\subsubsection{Tensor product of diagrams over a small category}

If $\mA$ is a cocomplete additive category, then for any small category $\mC$ there is defined a bifunctor {\it of a tensor product} that is additive in each of the arguments
$$
\otimes: \Ab^{\mC}\times \mA^{\mC^{op}} \to \mA,
$$
whose value is characterized by the existence
isomorphisms
$$
\mA(G\otimes F, A)\stackrel{\cong}\to \Ab^{\mC}(G, Hom(F,A)),
$$
natural in $G\in \Ab^{\mC}$, $F\in \mA^{\mC^{op}}$ and $A\in \mA$.

This bifunctor has the following property:
There is a natural isomorphism $\xi_c: \ZZ h^c\otimes F\stackrel{\cong}\to F(c)$.
Being natural means that for every morphism $\alpha: a\to b$ the diagram
\begin{equation}\label{nattrans}
\xymatrix{
\ZZ h^a\otimes F \ar[rr]^{\xi_a} && F(a)\\
\ZZ h^b\otimes F \ar[u]^{\ZZ h^{\alpha}\otimes 1_F} \ar[rr]^{\xi_b}
 && F(b)\ar[u]_{F(\alpha)}
}
\end{equation}
For fixed $F\in \mA^{\mC^{op}}$ the functor
 $(-)\otimes F$ commutes with colimits. For every
 $G\in \Ab^{\mC}$ the functor  $G\otimes (-)$ commutes with colimits
 \cite[Lemma 3.2]{X2008}.

\subsubsection{Constructing complexes for homology of categories}

\begin{proposition}
Let $\ZZ h^{(-)}: \mC^{op}\to \Ab^{\mC}$ be the composition of the functor $\ZZ(-): \Set\to \Ab$ and Yoneda embedding $h^{(-)}: \mC^{op}\to \Set^{\mC}$.
Then $\coLim_n^{\mC^{op}}\ZZ h^{(-)}=0$ for all $n\geq 1$, and $\coLim^{\mC^{op}}\ZZ h^{(-)}= \Delta_{\mC}\ZZ$.
\end{proposition}
{\sc Proof.}
The values $\coLim_n^{\mC^{op}}\ZZ h^{(-)}$ are isomorphic to the homology object of the complex $C_n(\mC^{op}, \ZZ h^{(-)})$ consisting of from objects and morphisms of the category $\mA^{\mC}$.
We need to prove the accuracy of the sequence
\begin{equation}\label{seq1}
0 \leftarrow \Delta_{\mC}\ZZ \leftarrow \bigoplus\limits_{c_0\in \mC}{\ZZ h^{c_0}}
\stackrel{d_1}\leftarrow \bigoplus\limits_{c_0\leftarrow c_1}{\ZZ h^{c_0}}
\stackrel{d_2}\leftarrow
\bigoplus\limits_{c_0\leftarrow c_1\leftarrow c_2}{\ZZ h^{c_0}}
\leftarrow \cdots
\end{equation}
in the category $\mA^{\mC}$. This sequence is exact if and only if for each $c\in \mC$ the sequence of values
$$
0 \leftarrow \ZZ \leftarrow \bigoplus\limits_{c_0\in \mC}{\ZZ\mC(c_0,c)}
\stackrel{(d_1)_c}\leftarrow \bigoplus\limits_{c_0\leftarrow c_1}{\ZZ\mC(c_0,c)}
\stackrel{(d_2)_c}\leftarrow
 \bigoplus\limits_{c_0\leftarrow c_1\leftarrow c_2}{\ZZ\mC(c_0,c)}
\leftarrow \cdots
$$
The resulting complex will consist of free Abelian groups with differentials defined on the bases of these groups by the formula
$(d_n)_c(c_n\to \cdots \to c_0\to c)=
\sum\limits_{i=0}^n (-1)^i (c_n\to \cdots \to \widehat{c_i}\to \cdots \to c_0\to c)$.
Here $\widehat{c_i}$ denotes the removal of the object $c_i$ from the sequence of morphisms, followed by the replacement, in the case of $0\leq i\leq n-1$, that goes out of $c_i$ and enters $c_i$ morphisms by their composition. In the case $i=n$, the morphism outgoing from $c_n$ is removed. The morphism $(d_0)_c: \bigoplus\limits_{c_0\in \mC}{\ZZ\mC(c_0,c)}\to \ZZ$ assigns to each element of the basis an element $1\in \ZZ$.
We now define the homomorphisms
$s_n: \bigoplus\limits_{c_0\leftarrow \cdots \leftarrow c_{n-1}}{\ZZ\mC(c_0,c)}
\to \bigoplus\limits_{c_0\leftarrow \cdots \leftarrow c_n}{\ZZ\mC(c_0,c)}$,
acting for  $n\geq 1$ on elements of bases by the formula
$s_n(c_{n-1}\to \cdots \to c_0\to c)=
(c_{n-1}\to \cdots \to c_0\to c\stackrel{1_c}\to c)$.
And we define the homomorphism $s_0: \ZZ\to \bigoplus\limits_{c_0}{\ZZ\mC(c_0,c)}$,
by setting $s_0(1)=(c\stackrel{1_c}\to c)$.
The equalities $(d_0)_c s_0=1_{\ZZ}$ and $(d_{n+1})_c s_{n+1}+s_n(d_n)_c=1$ hold.
This implies that $s$ will be a homotopy between the identity and zero morphisms of the complex into itself. Hence the sequence  (\ref{seq1}) is exact.
\hfill $\Box$

\begin{proposition}\label{projcolim}
For any projective resolution $P_*\to \Delta\ZZ$ in the category $\Ab^{\mC}$ there are 
isomorphims
$\coLim_n^{\mC^{op}}F\cong H_n(P_*\otimes F)$,
natural in $F\in \mA^{\mC^{op}}$.
\end{proposition}
{\sc Proof.}
Substituting the dual category for $\mC$, one can obtain the functor $C_*(\mC^{op},-): \mA^{\mC^{op}}\to Ch(\mA)$ into the category of chain complexes.
Since there is an isomorphism $\ZZ h^{(-)} \otimes F\stackrel{\cong}\to F(-)$, then we have the isomorphism $C_*(\mC^{op} , \ZZ h^{(-)}\otimes F)\stackrel{\cong}\to C_*(\mC^{op},F)$.
The complex $C_*(\mC^{op}, \ZZ h^{(-)}\otimes F)$ consists of diagrams
$$
\bigoplus\limits_{c_0\leftarrow \cdots \leftarrow c_n}\ZZ h^{c_0}\otimes F
\cong (\bigoplus\limits_{c_0\leftarrow \cdots \leftarrow c_n}\ZZ h^{c_0})\otimes F.
$$
This implies that the homology objects of the complex $C_*(\mC^{op},F)$ are isomorphic to the homology objects of the complex obtained by multiplying the projective resolution (\ref{seq1}) of the diagram $\Delta_{\mC}\ZZ$ by the diagram $F$:
$$
0 \leftarrow (\bigoplus\limits_{c_0\in \mC}{\ZZ h^{c_0}})\otimes F
\stackrel{d_1}\leftarrow (\bigoplus\limits_{c_0\leftarrow c_1}{\ZZ h^{c_0}})\otimes F
\stackrel{d_2}\leftarrow
(\bigoplus\limits_{c_0\leftarrow c_1\leftarrow c_2}{\ZZ h^{c_0}})\otimes F
\leftarrow \cdots.
$$
Any projective resolution $P_*$ of the diagram $\Delta_{\mC}\ZZ$ is homotopy equivalent to the resolution (\ref{seq1}).
This implies that the homology objects $H_n(P_*\otimes F)$ will be isomorphic to the homology objects of $C_*(\mC^{op},F)$, which, in turn, are isomorphic to the values of the left satellites $\coLim^ {\mC^{op}}_n F$ of the colimit functor.\hfill $\Box$

\begin{example}\label{hsimp}
Let $\mA= \Ab$ be the category of Abelian groups, and $F= \Delta_{\mC}\ZZ$.
According to the proposition \ref{projcolim}, for $\mD=\mC^{op}$, the homology groups of $C_*(\mD, \Delta_{\mD}\ZZ)$ will be isomorphic to the homology groups of a simplicial set (in the sense \cite[Appendix II, \S1.1]{gab1967}) of a nerve of the category $\mD$ . Therefore, we can denote them by $H^{simp}_n(\nr\mD)$ 
or $H^{simp}_n\mD$.
\end{example} 

\section{Homology of $\mD$-sets with coefficients in systems}

Let $\mD$ be a small category.
A $\mD$-set, or a presheaf of sets on $\mD$, is a functor $X: \mD^{op}\to \Set$.

We will consider the homology of $\mD$-sets with coefficients in contravariant systems of objects of the AB4-category $\mA$, where $\mD$ is an arbitrary small category.

\subsection{Kan extension along virtual discrete prefibrations}

A subcategory is called co-reflective if its embedding functor has a right adjoint functor \cite[\S IV.5]{mac1972}. Let $S: \mC \to \mD$ be a functor between small categories.
We will call it a virtual discrete prefibration if for every object $d\in \mD$ the category $d/S$ contains a discrete
co-reflective subcategory.
This will be true if and only if each connected component of the category $d/S$ has an initial object.

Let $K$ be a small category containing a co-reflective discrete subcategory. Let $init(K)$ denote this discrete subcategory as well as the set of its objects. For each $k\in \Ob K$ we denote by $i(k)$
the initial object of the connected component in $K$ containing $k$, and by $i_k: i(k)\to k$ we denote the unique morphism from $i(k)$ to $k$.
($i: K\to init(K)$ will be the right adjoint functor to the embedding $init(K)\subseteq K$).

\begin{proposition}\label{lemma36}\cite[Lemma 3.6]{X2008}
Let $\mA$ be a category with colimits, and let $S: \mC\to \mD$ be a virtual discrete
 prefibration. Then for any functor $F: \mC^{op}\to \mA$ the left Kan extension $Lan^{S^{op}}: \mD^{op}\to \mA$ is isomorphic to a functor taking on objects $ d\in \mD$ values $\coprod\limits_{\beta \in init(d/S)} FQ^{op}_d(\beta)$ and assigning each morphism $\alpha: d\to e$ of the category $ \mD$ is a morphism $(Lan^{S^{op}}F) (\alpha)$ of the category $\mA$ defined by the following diagram commutativity condition
$$
\xymatrix{
\coprod\limits_{\gamma \in init(e/S)} FQ^{op}_e(\gamma) \ar[rrr]^{(Lan^{S^{op}}F) (\alpha)}
&&& \coprod\limits_{\beta \in init(d/S)} FQ^{op}_d(\beta)\\
F Q^{op}_e(\gamma)= F Q^{op}_d(\gamma\circ\alpha) \ar[u]^(.4){in_{\gamma}}
\ar[rrr]^{F Q^{op}_d(i_{\gamma\circ\alpha})} &&& F Q^{op}_d(i(\gamma\circ\alpha))
\ar[u]_(.4){in_{i(\gamma\circ\alpha)}}
}
$$ 
Here $in_{\gamma}$ are morphisms of the coproduct cone.
\end{proposition}

\subsection{Complexes for computing homology of $\mD$-sets}

Let $\mD$ be a small category. Consider an arbitrary functor $X: \mD^{op}\to \Set$. Each element $x\in \coprod\limits_{d\in Ob\mD}X(d)$ corresponds to a natural transformation $h_{\dom x}\stackrel{\widetilde{x}}\to X$, where $\ dom x =d$ is an object of category $\mD$ for which $x\in X(d)$.
This natural transformation is constructed using the Yoneda isomorphism $X(d) \cong \Set^{\mD^{op}}(h_d, X)$.

Let $\mD/X$ be a category whose objects are the natural transformations $h_{\dom x}\stackrel{\widetilde{x}}\to X$ and the morphisms $\widetilde{x}_1\stackrel{\alpha}\to \widetilde{x}_2$ are defined by $\alpha:\dom x_1\to \dom x_2$ morphisms in $\mD$
 that make triangles commutative
$$
\xymatrix{
h_{\dom x_1} \ar[rd]_{h_{\alpha}} \ar[rr]^{\widetilde{x}_1} && X\\
& h_{\dom x_2} \ar[ru]_{\widetilde{x}_2}
}
$$
Here $h_{\alpha}$ is a natural transformation whose components $(h_{\alpha})_d: \mD(d, \dom x_1)\to \mD(d, \dom x_2)$ for all $d \in Ob\mD$ and $\beta\in \mD(d, \dom x_1)$ are defined by the formula $(h_{\alpha})_d(\beta)= \alpha\circ \beta$.

{\it Category of elements} of the presheaf $X: \mD^{op}\to \Set$ consists of the set of objects $x\in \coprod_{d\in \Ob\mD}X(d)$. Its morphisms between $ x_1, x_2\in \coprod_{d\in \Ob\mD}X(d)$ serve as triples $x_1 \xrightarrow{\alpha} x_2$ such that $\alpha\in \mD(\dom x_2, \ dom x_1)$ and $X(\alpha)(x_1)= x_2$.
We see that every morphism $x_1 \xrightarrow{\alpha} x_2$ of the element category is equal to
morphism $x_1 \xrightarrow{\alpha} X(\alpha)(x_1)$.

Category of elements of presheaf $X: \mD^{op}\to \Set$
is isomorphic to the category $(\mD/X)^{op}$, and we will denote these categories in the same way. Isomorphism is achieved by an anti-isomorphism that associates each morphism $x_1 \xrightarrow{\alpha} x_2$ of the element category with a morphism $\tilde{x_2} \xrightarrow{\alpha} \tilde{x_1}$ of the category $\mD/X$.

A contravariant system of objects of the category $\mA$ on $X\in \Set^{\mD^{op}}$ is an arbitrary functor $G: (\mD/X)^{op}\to \mA$.

Let $Q_X: \mD/X \to \mD$ be a functor assigning to each object $\tilde{x}: h_{\dom x}\in X$ an object $\dom x\in \Ob\mD$, and each morphism $\tilde{x}_1 \xrightarrow{\alpha} \tilde{x}_2$ is a morphism $\alpha: \dom x_1 \to \dom x_2$.

\begin{proposition}\label{contralan}
Let $\mA$ be an AB4-category and let
$G: (\mD/X)^{op}\to \mA$ is a contravariant system of objects in $\mA$. Then
\begin{enumerate}
\item\label{conten1} The functor $Lan^{Q^{op}_X}G: \mD^{op}\to \mA$ takes on the objects $d\in \mD$ the values $\bigoplus_{x\in X(d)} G(x)$, and the morphisms $Lan^{Q^{op}_X}G (\alpha)$, $\alpha\in \mD(b,a)$, are defined by the following diagram commutativity condition for every $x\in \coprod_{d\in \Ob\mD}X(d)$:
\begin{equation}\label{dgleftke}
\xymatrix{
\bigoplus\limits_{z\in X(a)} G(z) \ar[rr]^{(Lan^{Q^{op}_X}G)(\alpha)}
&& \bigoplus\limits_{z\in X(b)}G(z)
\\
G(x)\ar[u]^{in_x} \ar[rr]_{G(x\xrightarrow{\alpha}X(\alpha)(x))} && G(X(\alpha)(x) )\ar[u]_{in_{X(\alpha)(x)}}
}
\end{equation}
\item\label{conten2} There are natural isomorphisms $\coLim^{(\mD/X)^{op}}_n G \cong \coLim^{\mD^{op}}_n Lan^{Q^{op}_X}G $ for all $n\geq 0$.
\end{enumerate}
\end{proposition}
{\sc Proof.}
According to \cite[\S X.3, formula (10), p.240]{mac1972}, for any functor $S: \mC\to \mD$ between small categories and for an arbitrary functor $F: \mC\to \mA$ into the cocomplete category $\mA$, the value of the left
Kan extensions for the functor $F$ along the functor $S$ on the objects $d\in \mD$ are calculated by the formula $Lan^S F (d) = \coLim^{S/d}FQ_d$, where $Q_d: S/d \to \mC$ 
is the forget functor of the left fiber of $S$ over $d\in \mD$.
Recall that $Q_d(c\in \mC, \alpha: S(c)\to d)= c$ on objects, and $Q_d((c,\alpha)\xrightarrow{\gamma}(c',\alpha'))= (c\xrightarrow{\gamma} c')$ - on morphisms of the category $S/d$.

Moreover, if $\mA$ is an AB4-category, then the same is true for the left satellites of the Kan extension:
$$
Lan^S_n F (d) = \coLim^{S/d}_n FQ_d, \text{ for all } n\geq 0
$$
\cite[Appendix II, note 3.8]{gab1967}.

Consider the presheaf $X: \mD^{op}\to \Set$.
Let $d\in \Ob\mD$. Objects of category $d/Q_X$
consist of pairs of morphisms $(d \xrightarrow{\alpha} Q_X(\tilde{x}), h_{\dom x}\xrightarrow{\tilde{x}}X)$.
We will denote them as pairs $(\tilde{x}, \alpha)$.
The morphism $(\tilde{x_1}, \alpha_1)\xrightarrow{\alpha} (\tilde{x_2}, \alpha_2)$ 
is defined by the morphism $\alpha: \dom x_1 \to \dom x_2$, for which the triangle
 shown is commutative on the left in the following figure:

$$
\xymatrix{
& Q_X(\tilde{x_1}) \ar[dd]^{\alpha}\\
d \ar[ru]^{\alpha_1} \ar[rd]^{\alpha_2} &\\
& Q_X(\tilde{x_2})
} \qquad
\xymatrix{
h_{\dom x_1} \ar[dd]_{h_{\alpha}} \ar[rd]^{\tilde{x}_1} \\
 & X  \\
h_{\dom x_2} \ar[ru]_{\tilde{x}_2}
}
$$
(The morphism $\tilde{x}_1 \xrightarrow{\alpha} \tilde{x}_2$ is shown on the right.)
The existence of the morphism $(\tilde{x}_1, \alpha_1) \xrightarrow{\alpha} (\tilde{x}_2, \alpha_2)$ is equivalent to the commutativity of the diagram$$
\xymatrix{
& h_{\dom \tilde{x}_1} \ar[dd]^{\alpha} \ar[rd]^{\tilde{x}_1}\\
h_d \ar[ru]^{\alpha_1} \ar[rd]_{\alpha_2} && X\\
& h_{\dom\tilde{x}_2} \ar[ru]_{\tilde{x}_2} 
}
$$
Hence, if the objects $(\tilde{x}_1, \alpha_1)$ and $(\tilde{x}_2, \alpha_2)$ belong to the same connected component, then $\tilde{x}_1 h_{\alpha_1}= \tilde{x}_2 h_{\alpha_2}$.
This implies that the connected component of the category $d/Q_X$ containing the object $(\tilde{x}, \alpha)$ has the initial object $(\tilde{x}h_{\alpha}, 1_d)$. This object is equal to $\widetilde{X(\alpha)x}$ and hence there is a bijection between the set of connected components and the set $X(d)$.
The following diagram illustrates the initial object of the connected component containing the object $(\tilde{x}, 1_d)$:
$$
\xymatrix{
& h_{\dom \tilde{x}_1} \ar[dd]^{\alpha} \ar[rd]^{\tilde{x}h_{\alpha}}\\
h_d \ar[ru]^{1_d} \ar[rd]_{\alpha} && X\\
& h_{\dom\tilde{x}} \ar[ru]_{\tilde{x}} 
}
$$

Therefore, the objects $(\tilde{x}_1, \alpha_1)$ and $(\tilde{x}_2, \alpha_2)$ belong to the same connected component if and only if $\tilde{x}_1 h_{\alpha_1 }= \tilde{x}_2 h_{\alpha_2}$.

The set of initial objects of the connected components is equal to 
$$
init(d/Q_X)= \{(\tilde{x}, 1_d) | x\in X(d)\}.
$$
 Denote the initial object of the connected component containing the object $(\tilde{x}, \alpha)$ by
$i(\tilde{x}, \alpha)$, it is equal to $(1_d, \tilde{x}h_{\alpha})$. There is a single morphism from the initial object to the object $(\tilde{x}, \alpha)$. Denote it by $i_{(\tilde{x}, \alpha)}$.

To prove the commutativity of the diagram (\ref{dgleftke}), we use the proposition \ref{lemma36}.
Consider the functor $S= Q_X: \mD/X\to \mD$.
For each $d\in \mD$ the functor $Q_d: d/Q_X\to \mD/X$ assigns to each pair $(\tilde{z}, \alpha)$ an object $\tilde{z}: h_{\dom z}\to X$ from the category $\mD/X$.
For each morphism $\alpha: b\to a$ of the category $\mD$ we define the functor $(-)\circ \alpha: a/Q_X\to b/Q_X$ acting on objects as $(\tilde{z}, \alpha')\mapsto (\tilde{z}, \alpha'\alpha)$.

Let $\gamma\in init(a/Q_X)$. This means the existence of $z\in X(a)$ such that $\gamma = (\tilde{z}, 1_a)$.
The equality $FQ^{op}_a(\gamma)= FQ^{op}_b(\gamma\circ\alpha)= FQ^{op}_b(\tilde{z}, \alpha)$ is true. There is a unique morphism $i_{(\tilde{z}, \alpha)}: (\tilde{z}h_{\alpha}, 1_b)\to (\tilde{z}, \alpha)$.
For any functor $F: (\mD/X)^{op}\to \mA$ we obtain the morphism $FQ^{op}_b(i_{(\tilde{z}, \alpha)}): FQ^{op }_b(\tilde{z}, \alpha) \to FQ^{op}_b(\tilde{z}h_{\alpha}, 1_b)$. This morphism is equal to the morphism in the bottom row of the diagram from \ref{lemma36}, from which we arrive at the commutative diagram
$$
\xymatrix{
\bigoplus\limits_{z\in X(a)} G(\tilde{z}) \ar[rr]^{(Lan^{Q^{op}_X}G)(\alpha)} 
							&& \bigoplus\limits_{z\in X(b)}G(\tilde{z})\\
G(\tilde{x})\ar[u]^(.4){in_{\tilde{x}}}
 \ar[rr]_{G(\tilde{x}h_{\alpha}\xrightarrow{\alpha} \tilde{x})} &&
 G(\tilde{x}h_{\alpha})\ar[u]_(.4){in_{\tilde{x}h_{\alpha} }}
}
$$

Received
a commutative diagram is transformed into a commutative diagram (\ref{dgleftke}) by using the rule: if the functor $G$ is defined on the category of elements, then the morphisms $\tilde{x}_1 \xrightarrow{\alpha} \tilde{x}_2$ from $\mD/X$ are identified with morphisms $x_2 \xrightarrow{\alpha} x_1$ from the category of elements.
This rule leads to a description of the morphism $Lan^{Q^{op}_X}G(\alpha)$ in the proposition being proved.
 
It remains to prove the isomorphism $\coLim^{(\mD/X)^{op}}_n G \cong \coLim^{\mD^{op}}_n Lan^{Q^{op}_X}G $.
The shortest way is to use the Andr\'e spectral sequence \cite[Appendix 2, Theorem 3.6]{gab1967}
$$
	E^2_{p q}=	\coLim^{\mD^{op}}_p Lan^{Q^{op}_X}_q G \Rightarrow \coLim^{(\mD/X)^{op}}_{p+q}G
$$
Since each connected component of the category $d/Q_X$ has an initial object, then each connected component of the category $Q^{op}_X/d$ has a final object, which implies the precision of the functor $Lan^{Q^{op}_X}$. Therefore, this spectral sequence degenerates and gives isomorphisms
$\coLim^{\mD^{op}}_n Lan^{Q^{op}_X}G \cong \coLim^{(\mD/X)^{op}}_n G$ 
for all $n\geq 0$.\hfill$\Box$

\begin{example}
Under the conditions of Proposition \ref{contralan}, consider the case $\mA=\Ab$, the category of abelian groups.
The functor $Lan^{Q^{op}_X}G: \mD^{op}\to \Ab$ assigns to each object $a\in \mD$ the direct sum of Abelian groups $\bigoplus\limits_{z\in X( a)} G(z)$.
It assigns to each morphism $\alpha: b\to a$ a homomorphism
$$
(Lan^{Q^{op}_X}G) (\alpha): \bigoplus\limits_{z\in X(a)} G(z) \to \bigoplus\limits_{z\in X(b)} G(z),
$$
acting on elements of direct summands as
  $$
(x\in X(a), g\in G(x)) \mapsto (X(\alpha)(x), G(x\xrightarrow{\alpha} X(\alpha)(x))(g) \in G(X(\alpha)(x))).
$$
\end{example}

\begin{corollary}\label{covarran}
Let $\mA$ be an AB4*-category, and let $G: \mD/X\to \mA$ be a covariant system of objects in $\mA$. Then there are natural
 isomorphisms $\Lim^n_{\mD/X} G \cong \Lim^n_{\mD} Ran_{Q_X}G $ for all $n\geq 0$.
Moreover, the functor $Ran_{Q_X}G: \mD\to \mA$ takes on the objects $d\in \mD$
 the values $\prod_{x\in X(d)}G(x)$, and the morphisms $Ran_ {Q_X}G (\alpha)$,
  $\alpha\in \mD(b,a)$, are defined by the following diagram commutativity condition 
for all $x\in \coprod_{d\in \Ob\mD}X(d)$\begin{equation}\label{dgrightke}:
\xymatrix{
\prod\limits_{z\in X(b)} G(z) \ar[d]_(.55){pr_{X(\alpha)x}} 
\ar[rr]^{(Ran_{Q_X}G)(\alpha)} 
						&& \prod\limits_{z\in X(a)}G(z) \ar[d]^(.55){pr_x} \\
G(X(\alpha)x) \ar[rr]_{G(x\xrightarrow{\alpha}X(\alpha)(x))} && G(x)
}
\end{equation}
\end{corollary}
{\sc Proof.} In Proposition \ref{contralan} we substitute $\mA^{op}$
instead of $\mA$ and $G^{op}$ instead of $G$. The suggestion will give a diagram and formulas in $\mA^{op}$. The transformation of this diagram and formulas into a diagram and a formula into $\mA$ will lead to the required corollary.
\hfill$\Box$

\begin{example}
Let us describe the formula for the right Kahn extension under the conditions of Corollary \ref{covarran} in the case $\mA=\Ab$.
The elements of the product $\prod_{z\in X(a)} G(z)$ will be considered as functions $\varphi: X(a)\to \cup_{z\in X(a)}G(z)$ satisfying condition $\varphi(z)\in G(z)$. The projection acts like $pr_x(\varphi)=x$.
We obtain for $\psi\in \prod_{z\in X(b)} G(z)$:
$$
(Ran_{Q_X}G)(\alpha)(\psi)(x)= G(x \xrightarrow{\alpha} X(\alpha)x)\psi(X(\alpha)x)
$$
\end{example}

\begin{remark}\label{fibrgroth}

The functor $Q_X: \mD/X\to \mD$ is a discrete Grothendieck fibration in which the fiber
 $Q^{-1}_X(d)$, $d\in \Ob\mD$, is equal to the image of the embedding $X(d)$ 
 in $d/Q_X$ assigning to each $x\in X(d)$ an object of the category $d/Q_X$ equal to 
 $(d\xrightarrow{1_d}Q_X(h_d\xrightarrow{\tilde{x}}X ))$.

This property of the functor $Q_X$ also allows us to prove
Proposition \ref{contralan}.
There are other problems that can be solved using discrete Grothendieck fibrations \cite{gal2012}, \cite{gal2013}, \cite{gal2021}.

To construct a commutative diagram (\ref{dgleftke}), it was easier for us to use the property of the existence of initial objects in the components of the right fibers of the functor $Q_X$.
 \end{remark}

\subsection{Direct image homology for a contravariant system of objects}

Let $X,Y: \mD^{op}\to \Set$ be set diagrams. For an arbitrary natural transformation 
$f: X\to Y$, denote by $\mD/f: \mD/X\to \mD/Y$ the functor that assigns to each object
 $h_d\stackrel{\widetilde{x}}\to X $ of category $\mD/X$ an object of category 
 $\mD/Y$ equal to composition $h_d\stackrel{f\circ\widetilde{x}}\to Y$. 
 This functor assigns to each morphism $( \alpha, \widetilde{x}_1, \widetilde{x}_2)$ a morphism $(\alpha: f\circ\widetilde{x}_1\to f\circ\widetilde{x} _2)$.
  Objects of the category $\widetilde{y}/ (\mD/f)$ are specified as pairs
   $(\widetilde{x}\in \mD/X, \alpha\in \mD(\dom y, \dom x))$ , for which the diagrams
$$
 \xymatrix{
 h_{\dom x}\ar[d]_{\widetilde{x}}
 & h_{\dom y}\ar[l]_{h_{\alpha}}\ar[d]_{\widetilde{y}}\\
 X \ar[r]_{f} & Y
 }
 $$
The morphisms $(\widetilde{x_1}, \alpha_1) \to (\widetilde{x_2}, \alpha_2)$ are given by the morphisms $\gamma: \dom x_1 \to \dom x_2$ making the following diagrams commutative
$$
\xymatrix{
h_{\dom{x_1}}\ar[dd]_{h_{\gamma}} \ar[rd]_{\widetilde{x}_1} \\
 &  X \ar[r]^f & Y & h_{\dom y}\ar[l]_{\widetilde{y}}
 \ar@(u,l)[lllu]_{h_{\alpha_1}} \ar@(d,l)[llld]^{h_{\alpha_2}}\\
h_{\dom{x_2}}\ar[ru]^{\widetilde{x}_2}
}
$$
Hence, if there exists a morphism $(\widetilde{x_1}, \alpha_1) \to (\widetilde{x_2}, \alpha_2)$, then $\widetilde{x_1}\circ h_{\alpha_1}= \widetilde{x_2} \circ h_{\alpha_2}$.
This implies the following assertion.
 \begin{lemma}\label{compcomma}
 For every $\widetilde{y}\in Ob(\mD/Y)$ the connected component of the comma category
  $\widetilde{y}/(\mD/f)$ containing the object $(\widetilde{x}, \alpha )$, has the initial object $i(\widetilde{x}, \alpha)= (\widetilde{x}\circ h_{\alpha},1_{\dom{y}})$. 
  The set of initial objects $init(\widetilde{y}/(\mD/f))$ of connected components 
   is equal to 
  $\{(\widetilde{x}, 1_{\dom x}) | \widetilde{x}\in \mD/X~ \& ~f\circ \widetilde{x}= \widetilde{y}\}$.
The unique morphism $i_{(\widetilde{x}, \alpha)}:
i(\widetilde{x}, \alpha)\to (\widetilde{x}, \alpha)$ is given by the morphism $\alpha$.
The forgetting comma-category functor $Q_{\widetilde{y}}: \widetilde{y}/(\mD/f)\to \mD/X$ acts on objects as $Q_{\widetilde{y}}(\widetilde {x}, \alpha)= \widetilde{x}$, and on morphisms as 
 $Q_{\widetilde{y}}((\widetilde{x}_1, \alpha_1)
\stackrel{\gamma}\to (\widetilde{x}_2, \alpha_2))= (\widetilde{x}_1
\stackrel{\gamma}\to \widetilde{x}_2)$.
\end{lemma}

The direct image of a contravariant system $F: (\mD/X)^{op}\to \mA$ of objects of cocomplete category $\mA$ on $X$ is called a contravariant system,
equal to the left Kan extension $Lan^{(\mD/f)^{op}}F$ \cite{mac1972}.

\begin{proposition}\label{landiscr}
For an arbitrary natural transformation of set diagrams $f: X\to Y$ and a contravariant system $F: (\mD/X)^{op}\to \mA$ in the cocomplete category $\mA$, the direct image is isomorphic to the contravariant system $f_* F$, which is defined on objects as $f_*F(\widetilde{y})= \coprod\limits_{x\in f^{-1}_{\dom y}(y)} F(\widetilde{x} )$,
and on the morphisms $\widetilde{y} \stackrel{\beta}\to \widetilde{z}$ in such away that the following diagrams are commutative:
$$
\xymatrix{
\coprod\limits_{x\in f^{-1}_{\dom z}(z)} F(\widetilde{x})
\ar[rr]^{(f_*F)(\widetilde{y}\stackrel{\beta}\to\widetilde{z})}
&& \coprod\limits_{x\in f^{-1}_{\dom y}(y)} F(\widetilde{x})\\
 F(\widetilde{x}) \ar[u]^(.35){in_{\widetilde{x}}}
\ar[rr]_{F(\widetilde{x}\circ
 h_{\beta}\stackrel{\beta}\to\widetilde{x})}
 &&  F(\widetilde{x}\circ h_{\beta})\ar[u]_(.35){in_{\widetilde{x}\circ h_{\beta}}}
}
$$
\end{proposition}
{\sc Proof.} The set $init(\widetilde{y}/(\mD/f))$ consists of pairs 
$(\widetilde{x}, 1_{\dom y})$ satisfying $f\circ \widetilde{x}= \widetilde{y}$ 
which is
equivalent to $x\in f^{-1}_{\dom y}(y)$. Hence Lemma \ref{compcomma} 
on the structure of the comma-category together with Proposition \ref{lemma36} 
lead to the direct image isomorphism and a functor that takes on $\widetilde{y}\in \mD/Y$ 
the values 
$$
\coprod\limits_{x\in f^{-1}_{\dom y}(y)} FQ^{op}_{\widetilde{y}}(\widetilde{x},
1_{\dom y})= \coprod\limits_{x\in f^{-1}_{\dom y}(y )} F(\widetilde{x}).
$$
The diagram from Proposition \ref{lemma36} for defining the action of the left Kan extension
of the functor $F$ on morphisms leads to a diagram defining the functor $f_*F$. 
This implies that the constructed functor $f_* F$ will be isomorphic to the functor
$Lan^{(\mD/f)^{op}}F$.
\hfill $\Box$

\begin{theorem}\label{lanhomol}
Let $\mA$ be an AB4-category, $\mD$ a small category,
 $f: X\to Y$ a morphism between $X, Y\in \Set^{\mD^{op} }$. Then for an arbitrary functor $F: (\mD/X)^{op}\to \mA$ there are natural isomorphisms
 $\coLim^{(\mD/X)^{op}}_n F \stackrel{\cong}\to \coLim^{(\mD/Y)^{op}}_n f_*F$.
\end{theorem}
{\sc Proof.} Consider the Andre spectral sequence described in
 \cite[Appendix II, Theorem 3.6]{gab1967} applied to the functor
  $(\mD/f)^{op}: (\mD/X)^{op }\to (\mD/Y)^{op}$,
$$
E^2_{p,q}= \coLim^{(\mD/Y)^{op}}_p (Lan^{(\mD/f)^{op}}_q F) \Rightarrow \coLim^{( \mD/X)^{op}}_{p+q} F.
$$
Using \cite[Appendix II, Note 3.8]{gab1967}, we get $(Lan^{(\mD/f)^{op}}_q F)(\widetilde{y})\cong
\coLim^{(\widetilde{y}/(\mD/f))^{op}}_q (F\circ Q^{op}_{\widetilde{y}})$.
Each connected component of the category ${(\widetilde{y}/(\mD/f))^{op}}$ has a terminal object, which implies that the colimit functor over this category is exact. We get $(Lan^{(\mD/f)^{op}}_q F)(\widetilde{y})=0$, for $q>0$.
Hence the spectral sequence degenerates into the string $E^2_{p,0}$.
Hence, $\coLim^{(\mD/X)^{op}}_n F \stackrel{\cong}\to \coLim^{(\mD/Y)^{op}}_n (Lan^{ (\mD/f)^{op}} F)$.
According to the proposition \ref{landiscr}, there is an isomorphism $Lan^{(\mD/f)^{op}} F \cong f_*F$.
We arrive at the isomorphism $\coLim^{(\mD/X)^{op}}_n F \stackrel{\cong}\to \coLim^{(\mD/Y)^{op}}_n f_*F$.
\hfill $\Box$

\subsection{Criterion for isomorphic homology of $\mD$-sets}

Let $f: X\to Y$ be a morphism of $\mD$-sets.
The functor $\mD/f: \mD/X\to \mD/Y$ maps each object $\widetilde{x}: h_d\to X$ to a composition
$f\circ\widetilde{x}: h_d\to Y$.
Let $f^*: \mA^{(\mD/Y)^{op}} \to \mA^{(\mD/X)^{op}}$ be a functor that associates each functor with $F: ( \mD/Y)^{op}\to \mA$ composition $F\circ(\mD/f)^{op}$.
Since the functors $\coLim^{(\mD/Y)^{op}}_n: \mA^{(\mD/Y)^{op}}\to \mA$ are satellites of the colimit functor, there exists a unique sequence of canonical morphisms constituting the $\partial$-functor morphism
$$
\coLim^{(\mD/X)^{op}}_n F\circ(\mD/f)^{op}\to
\coLim^{(\mD/Y)^{op}}_n F.
$$
The question arises: when will these morphisms be isomorphisms?
To give an answer, let's use Oberst's theorem \cite[Theorem 2.3]{obe1968}.

Denote by $\overleftarrow{f}({y})$ the $\mD$-set defined by a pullback in the category of $\mD$-sets
\begin{equation}\label{invimage}
\xymatrix{
X \ar[r]^f & Y \\
\overleftarrow{f}({y}) \ar[r] \ar[u]_{f_y}
& h^{\mD}_d\ar[u]_{\widetilde{y}}
}
\end{equation}

This $\mD$-set $\overleftarrow{f}(y)$ is called
{\it inverse fiber} over the element $y\in Y(d)$.
It is easy to see that there is an isomorphism of categories
$(\Box/f)/\widetilde{y}\cong \mD/\overleftarrow{f}(y)$.
Applying Oberst's theorem, we arrive at the following statement. Groups $H^{simp}_n\mC$
are defined in Example \ref{hsimp}.

\begin{proposition}\label{critiso}
Let $f: X\to Y$ be a morphism of $\mD$-sets.
Then the following properties of the morphism $f$ are equivalent:
\begin{enumerate}
\item For every $y\in Y$ the category $\mD/\overleftarrow{f}(y)$ is connected and the homology groups
 $H^{simp}_n(\mD/\overleftarrow{f}(y)))$ are equal to zero for all $n> 0$.
\item The canonical homomorphisms $\coLim^{(\mD/X)^{op}}_n f^{*}F\to \coLim^{(\mD/Y)^{op}}_nF$ are isomorphisms for every functor $F: (\mD/Y)^{op}\to \Ab$.
\item The canonical morphisms $\coLim^{(\mD/X)^{op}}_n f^{*}F\to \coLim^{(\mD/Y)^{op}}_nF$ are isomorphisms for any functor $F: (\mD/Y)^{op}\to \mA$ into an arbitrary  AB4-category $\mA$.
\end{enumerate}
\end{proposition}

\subsection{Spectral sequence for the colimit homology of $\mD$-sets}

Let us find conditions for the diagram of $\mD$-sets under which there exists a spectral sequence converging to the homology of its colimit.

We will consider the spectral sequences of the first quarter, in the sense of \cite{mac1963}.
The following assertion is dual to \cite[Corollary 2.4]{X1989} obtained with \cite[Theorem 2.1]{X1989}. Instead of the AB4-category, 
the category $\mA^{op}$ with exact coproducts must be substituted.

\begin{proposition}\label{covering}
Let $J$ be a small category, and let $\{X^i\}_{i\in J}$ be a diagram of $\mD$-sets such that
\begin{equation}\label{c3}
\coLim_q^J\{\ZZ(X^i(d))\}_{i\in J}= 0, \mbox{ for any } d\in \Ob\mD \mbox{ and }
q>0.
\end{equation}
Let $\lambda_i: X^i\rightarrow \coLim^J\{X^i\}_{i\in J}$ be the cone of morphisms for the colimit of $\mD$-sets.
Then for any AB4-category $\mA$ and any contravariant system $F: (\mD/\coLim^J\{X^i\})^{op}\to \mA$ there exists a spectral sequence of the first quarter
$$
E^2_{p,q}= \coLim_p^{J}\{\coLim^{(\mD/X^i)^{op}}_q \lambda^*_i F\}_{i\in J} \Rightarrow
\{\coLim^{(\mD/\coLim^J\{X^i\})^{op}}_{p+q}F \}.
$$
\end{proposition}

\subsection{Spectral sequence of $\mD$-set morphism} 

Let $f: X\rightarrow Y$ be a morphism of $\mD$-sets. For each $\sigma\in {\mD}/Y$ the Cartesian square (\ref{invimage}) defines an inverse fiber over $\sigma$ and a morphism $f_{\sigma}: \overleftarrow{f}(\sigma) \rightarrow X$.
Inverse fibers over elements from $Y$ form a diagram of $\mD$-sets $\{\overleftarrow{f}(\sigma)\}_{\sigma\in {\mD}/Y}$ whose colimit is isomorphic to $X $. Applying the general result on the spectral sequence of the morphism \cite[Theorem 4.1]{X1991}, where the category $\mA$ should be replaced by the category $\mA^{op}$, we obtain
the following statement:

\begin{proposition}\label{spmor}
Let $f: X\rightarrow Y$ be a morphism of $\mD$-sets, $F$ be a contravariant system in the 
AB4-category $\mA$ on $X$. Then there is a spectral sequence of the first quarter
$$
E^2_{p,q}= \coLim_p^{\mD/Y}\{\coLim^{(\mD/\overleftarrow{f}(\sigma))^{op}}_q 
f^*_{\sigma} F\}_{\sigma\in \mD/Y} \Rightarrow
\{\coLim^{(\mD/X)^{op}}_{p+q}F \}.
$$
\end{proposition}
\begin{remark}\label{spmorinv}
If some diagram $G: \mC\to \mA$ consists of isomorphisms, then by inverting its morphisms we obtain a diagram on $\mC^{op}$. Denote it by $G^{-1}: \mC^{op}\to \mA$. In this case there are isomorphisms $\coLim^{\mC}_n G \cong \coLim^{\mC^{op}}_n G^{-1}$ \cite[Appendix II, Proposition 4.4]{gab1967}.
In particular, if the diagram
$$
\{\coLim^{(\mD/\overleftarrow{f}(\sigma))^{op}}_q
f^*_{\sigma} F\}_{\sigma\in \mD/Y}
$$
consists of isomorphisms, then by inverting these isomorphisms we obtain the diagram
$
\{\coLim^{(\mD/\overleftarrow{f}(\sigma))^{op}}_q
f^*_{\sigma} F\}_{\sigma\in \mD/Y}^{-1}
$.
Proposition \ref{spmor} leads to a spectral sequence 
$$
E^2_{p,q}\cong \coLim_p^{(\mD/Y)^{op}}\{\coLim^{(\mD/\overleftarrow{f}(\sigma))^{op}}_q 
f^*_{\sigma} F\}_{\sigma\in \mD/Y}^{-1} \Rightarrow
\{\coLim^{(\mD/X)^{op}}_{p+q}F \},
$$
 connecting the homology objects of $\mD$-sets $X$ and $Y$.
\end{remark}

\section{Homology of the cube category}

Consider the cube category and study the left satellites of the colimit functor 
$\coLim_n^{\Box^{op}}: \mA^{\Box^{op}}\to \mA$
on cubical objects in the AB4-category $\mA$.

\subsection{Cube category}
For an arbitrary $n\in \NN$, we will consider a partially ordered set $\II^n=\{0,1\}^n$ 
equal to the Cartesian power of a linearly ordered set $\II=\{0,1\}$. For $n=0$ the set 
$\II^0$ consists of the only element $\emptyset$. A partially ordered set $\II^n$ is called 
an $n$-dimensional {\it cube}.

The objects of the {\it category of cubes} $\Box$ are the cubes $\II^0$, $\II^1$, $\II^2$, \ldots.

Morphisms of the cube category are defined as nondecreasing mappings of these posets that can be decomposed into a composition of mappings of the form $\delta_i^{k,\varepsilon}: \II^{k-1}\rightarrow \II^k$ and $ \sigma_i^k: \II^k\to \II^{k-1}$
 defined for $k\geq 1$, $1\leq i\leq k$, $\varepsilon\in \{0,1\}$, and taking values

\begin{gather}
\delta_i^{k,\varepsilon}(x_1, \dots, x_{k-1})=
(x_1, \dots, x_{i-1}, \varepsilon, x_i, \dots, x_{k-1}),\\
\sigma_i^k(x_1, \cdots, x_k)=(x_1, \cdots, x_{i-1}, x_{i+1}, \cdots, x_k).
\end{gather}
In particular, $\delta_1^{1,0}(\emptyset)=0$, $\delta_1^{1,1}(\emptyset)=1$,
$\sigma^1_1(x_1)=\emptyset$, for all $x_1\in \II$.

The category of cubes is described in \cite{jar2006}, \cite{bro2011}.
In \cite{bro2011}, its objects are the cubes $[0,1]^n$.
\cite{jar2006} presents commutative diagrams illustrating relations that can be used 
to define the category $\Box$.

There are relations
\begin{gather}
\label{rel1}\delta_j^{n,\beta}\delta_i^{n-1,\alpha}=
\delta_i^{n,\alpha}\delta_{j-1}^{n-1,\beta}\quad  (1\leq i< j \leq n,
\alpha\in \II, \beta\in \II);\\
\label{rel2}
\sigma^{n-1}_j \sigma^n_i= \sigma^{n-1}_i\sigma^n_{j+1}\quad (1\leq i\leq j\leq n-1,
n\geq 2);\\
\label{rel3}
\sigma^{n+1}_j\delta^{n+1,\alpha}_i=\left\{
\begin{array}{lc}
\delta^{n,\alpha}_i\sigma^{n}_{j-1}, & (1\leq i<j\leq n+1, \alpha\in \{0,1\}), \\
\delta^{n,\alpha}_{i-1}\sigma^{n}_j, & (1\leq j<i\leq n+1, \alpha\in \{0,1\}),\\
1_{\II^n}, & i=j.
\end{array}
\right.
\end{gather}

According to \cite[Lemma 4.1]{gra2003}, every morphism of this category
 $f: \II^k\to \II^n$ admits the canonical composition\begin{equation}\label{grandis}
f= \delta^{n,\varepsilon_1}_{j_1}\cdots \delta^{k-r+1,\varepsilon_s}_{j_s}
\sigma^{k-r+1}_{i_1}\cdots \sigma^k_{i_r},
\quad
\begin{array}{l}
 1\leq i_1< \cdots <i_r\leq k,\\
 n\geq j_1>\cdots >j_s\geq 1,\\
 k-r=n-s\geq 0.
\end{array}
\end{equation}

\subsection{Cubical sets}
Let $\mA$ be an arbitrary category. {\it A cubical object} in the category $\mA$ is a functor $X: \Box^{op}\to \mA$.
{\it A cubical set} is a functor $X: \Box^{op}\to \Set$.

The category $\Box$ can be defined using a {\it graph with relations}, in the sense of \cite{bor1994}. In our case, this graph has vertices $\II^n$, where $n$ runs through all non-negative integers. Its edges are {\it boundary morphisms} $\II^{n-1} \stackrel{\delta^{n,\varepsilon}_i}\to \II^{n}$ and {\it degeneration morphisms} $\II^n\stackrel{\sigma^n_i}\to \II^{n-1}$, where the indices range over the values $n\geq 1$, $1\leq i\leq n$, $\varepsilon\in 
\{ 0.1\}$.
The relations are given above (\ref{rel1})-(\ref{rel3}).

The cubical object $X: \Box^{op}\to \mA$ can be given as
a set of objects $X_n=X(\II^n)$, $n\geq 0$, and a set of morphisms
$\partial^{n,\varepsilon}_i=X(\delta^{n,\varepsilon}_i): X_n\to X_{n-1}$,
$\frak{s}^n_i=X(\sigma^n_i): X_{n-1}\to X_n$, as well as a set of commutative diagrams,
corresponding to the relations (\ref{rel1})-(\ref{rel3}).

Morphisms $\partial^{n,\varepsilon}_i$ are called {\it boundary operators},
and $\frak{s}^n_i$ are {\it degeneration operators}.

\subsection{Non-degenerate cubes of the standard cube}

Let $(X_n, \partial^{n,\varepsilon}_i, \frak{s}^n_i)$ be a cubical set. A cube $x\in X_n$ is 
called {\it degenerate} if there are $y\in X_{n-1}$ and $i\in\{1, \ldots, n\}$ such that $\frak{ s}^n_i(y)=x$.
Otherwise, it is called {\it non-degenerate}.

Let $h_{\II^n}: \Box^{op}\to \Set$ be a morphism functor; it acts on objects as $h_{\II^n}(\II^k)= \Box(\II^k,\II^n)$, and to each morphism $f: \II^m\to \II^k$ it associates 
the natural transformation $\Box(f,\II^n): \Box(\II^ k,\II^n)\to \Box(\II^m,\II^n)$ 
assigning each $g\in \Box(\II^k,\II^n)$ an element $gf\in \Box(\II^m,\II^n)$.
The functor $h_{\II^n}$ is a cubical set and is called the {\it standard cube}. A standard cube can be defined as a triple $(X_k, \partial^{k,\varepsilon}_i, \frak{s}^k_i)$ consisting of a sequence of sets $X_k= \Box(\II^k,\II^n )$ in place with mappings $\partial^{k,\varepsilon}_i(f)= f\delta^{k,\varepsilon}_i$ and $\frak{s}^k_i(g)= g\sigma^k_i $.

According to the definition of a degenerate cube, the cube $f\in h_{\II^n} \II^k)$ is degenerate if and only if there exist
$g: \II^{k-1}\to \II^n$ and $i$ from the interval $1\leq i\leq k$ such that $f= g\sigma^k_i$.

\begin{proposition}\label{nondeg}\cite[Proposition 3]{X2019}
For every $n\geq 0$, the cube $f\in h_{\II^n}(\II^k)$ is non-degenerate if and only if the map $f: \II^k\to \II^n$ is an injection.
\end{proposition}

\subsection{Construction of a projective resolution}

Consider a sequence of morphisms in the category $\Ab^{\Box}$
\begin{equation}\label{complex}
0 \stackrel{d_0}\leftarrow \ZZ h^{\II^0} \stackrel{d_1}\leftarrow
\ZZ h^{\II^1} \stackrel{d_2}\leftarrow \ZZ h^{\II^2} \stackrel{d_3}\leftarrow
\ZZ h^{\II^3}
\leftarrow \cdots
\end{equation}
consisting of cocubical abelian groups and natural transformations 
$$
d_k= \sum\limits^k_{i=1}(-1)^i(\partial^{k,0}_i-\partial^{k,1}_i).
$$
Here $\partial^{k,\varepsilon}_i: \ZZ h^{\II^k}\to \ZZ h^{\II^{k-1}}$ are natural transformations whose components $(\partial^{k,\varepsilon}_i)_{\II^n}:\ZZ\Box(\II^k,\II^n) \to \ZZ\Box(\II^{k-1},\II^n)$ for objects $\II^n\in \Box$ are defined on elements of the basis $f\in \Box(\II^k,\II^n)$ by the formula $(\partial^{k,\varepsilon}_i)_{\II^n}(f)= f\delta^{k,\varepsilon}_i$.
From the relations $\partial^{k-1,\alpha}_i\partial^{k,\beta}_j= \partial^{k-1,\beta}_{j-1}\partial^{k,\alpha}_i$ resulting from the equality (\ref{rel1}) will be followed by $d_k d_{k+1}=0$ for all $k\geq 0$.
Therefore, the sequence of morphisms (\ref{complex}) is a chain complex.

For each $n\geq 0$, consider the subset $D_k(\II^n) \subseteq \Box(\II^k, \II^n)$ consisting of degenerate cubes of the cubical set $h_{\II^n}$ . The morphisms $\Box(\II^k, \delta^{n,\varepsilon}_i)$ and $\Box(\II^k, \sigma^{n}_i)$ carry degenerate cubes to degenerate ones, whence $D_k $ will be a subfunctor of the functor $h^{\II^k}$. Consider an embedding of functors $\ZZ D_k\subseteq \ZZ h^{\II^k}$, for an arbitrary $k\geq 0$. The cokernel of this embedding will be the functor $\ZZ h^{\II^k}/\ZZ D_k$ taking values on objects equal to the quotient groups  $\ZZ h^{\II^k}(\II^n)/\ZZ D_k(\II^n)$. The components of the projection 
$\ZZ h^{\II^k}\to \ZZ h^{\II^k}/\ZZ D_k$ are equal to the 
canonical projections onto quotient groups.

There is an exact sequence in $\Ab^{\Box}$
\begin{equation}\label{factor}
0 \leftarrow \ZZ h^{\II^k}/\ZZ D_k \stackrel{pr}\leftarrow \ZZ h^{\II^k} \stackrel{\supseteq}\leftarrow \ZZ D_k \leftarrow 0.
\end{equation}

\begin{lemma}\label{proj}\cite[Lemma 4]{X2019}
The functor $\ZZ h^{\II^k}/\ZZ D_k$ is a projective object of the category $\Ab^{\Box}$, for all integers $k\geq 0$.
\end{lemma}
{\sc Proof.} In \cite{X2019} a proof is given that contains gaps in the assertion that the section of $pr_{\II^n}$ mappings is natural.
In order to fix this, we construct a natural transformation $r: \ZZ h^{\II^k}\to \ZZ D_k$, inverse from the left to the embedding $\ZZ D_k \subseteq \ZZ h^{\II^k} $.
By Yoneda's lemma, in order to construct $r$, it suffices to specify an element $z= r_{\II^k}(1_{\II^k})\in \ZZ D_k(\II^k)$. And then the natural transformation $r$ will have components $r_{\II^n}(\alpha)=\alpha\circ z$, for all $n\geq 0$ and $\alpha\in\Box(\II^ k, \II^n)$.
Using the idea of Eilenberg and MacLane used in \cite[Prop. 7.2]{eil1953} in proving the representability of normalized groups of singular cubical homology, we set
$$
z = 1 - (1- \delta^{k,0}_1\sigma_1)(1- \delta^{k,0}_2\sigma_2)
\cdots(1- \delta^{k,0}_k\sigma_k).
$$
This $z$ is equal to the linear combination of products $ \delta^{k,0}_{s_1}\sigma_{s_1}\delta^{k,0}_{s_2}\sigma_{s_2}\cdots
\delta^{k,0}_{s_m}\sigma_{s_m}$ for some $1\leq m\leq k$ and $s_1 < \ldots < s_m$.
Each of these products will not be a monomorphism, since multiplying this product from the right by $(1-\delta^{k,0}_{s_m} \sigma_{s_m})$ gives $0$. Hence $z$ is equal to a linear combination of degenerate morphisms, and $z\in \ZZ D_k(\II^k)$.
Therefore, for any $n\geq 0$ and $\alpha\in \ZZ h^{\II^k}(\II^n)$, the element $r_{\II^n}(\alpha)$ belongs to $\ ZZ D_k(\II^n)$.

For any $1\leq i < j \leq k$ the operations $\delta^{k,0}_i\sigma^k_i$
and $\delta^{k,0}_j\sigma^k_j$ commute.
If $\alpha\in D^k(\II^n)$ then $\alpha= x\sigma_i$, for some $1\leq i\leq k$ and $x\in h^{\II^k} (\II^n)$.
By virtue of the permutability property noted,
$z= 1 - (1- \delta^{k,0}_i\sigma_i)y$ for some $y\in \Box(\II^k, \II^k)$. It entails
$$
x\sigma^k_i\circ z= x\sigma^k_i(1 - (1- \delta^{k,0}_i\sigma^k_i)y)=
x\sigma^k_i,
$$
and hence $r_{\II^n}(\alpha)= \alpha$ for all $\alpha\in \ZZ D_k(\II^n)$.

Hence $r$ is a retraction of $\ZZ h^{\II^k}$ onto $\ZZ D_k$, and
the short exact sequence (\ref{factor}) is split. Since $\ZZ h^{\II^k}$ is a projective object in $\Ab^{\Box}$, then $\ZZ h^{\II^k}/\ZZ D_k$ is projective.

The section $s: \ZZ h^{\II^k}/ \ZZ D_k \to \ZZ h^{\II^k}$ of the natural transformation $pr$ is determined from the retraction $r: \ZZ h^{\II^ k}\to \ZZ D_k$ standard \cite[\S I.4]{mac1963}, according to the formula
$
s_{\II^n}(\alpha+ \ZZ D_k(\II^n))= \alpha -\alpha z.
$
\hfill$\Box$

\begin{remark}
Relations between morphisms of the cube category involved in the proof that the functor 
$\ZZ h^{\II^k}/\ZZ D_k$ is projective were used by A. \`Swi\c{a}tec \cite{swi1981} for study of cubical objects in the abelian category and  I. Pachkoria \cite{pat2012} to study
pseudo-cubical objects of an idempotently complete preadditive category.
\end{remark}

\begin{remark}
We have obtained formulas for projections and embeddings that are natural in $\II^n$\begin{equation}\label{cosplit}
\ZZ D_k(\II^n) 
\begin{smallmatrix}
\xleftarrow{r_{\II^n}}\\
\subseteq
\end{smallmatrix}
\ZZ h^{\II^k}(\II^n)
\begin{smallmatrix}
\xrightarrow{\pi_n}\\
\leftarrowtail\\
s_n
\end{smallmatrix} 
\ZZ h^{\II^k}/ \ZZ D_k(\II^n)
\end{equation}
where
$r_{\II^n}(\alpha)=\alpha z$, $\pi_n(\alpha)= \alpha-\alpha z+ \ZZ D_k(\II^n)= 
\alpha+ \ZZ D_k(\II^n)$,
$s_n(\alpha+\ZZ D_k(\II^n))= \alpha-\alpha z$, $\alpha\in \ZZ h^{\II^k}(\II^n)$.
We have proved that $r_{\II^n}(\alpha)= \alpha$ is true for all $\alpha\in \ZZ D_k(\II^n)$.

Let's check the correctness of $s_n$ mapping:
\begin{multline*} 
\alpha_1+ \ZZ D_k(\II^n)= \alpha_2+ \ZZ D_k(\II^n)
\Rightarrow \alpha_1- \alpha_2\in \ZZ D_k(\II^n)\\
\Rightarrow (\alpha_1-\alpha_2)z= \alpha_1-\alpha_2 
\Rightarrow  \alpha_1-\alpha_1 z = \alpha_2-\alpha_2 z.
\end{multline*} 

Projections and embeddings (\ref{cosplit}) lead to a decomposition of the free cocubical abelian group $\ZZ h^{\II^k}$ into a direct sum of $\ZZ D_k$ and $\ZZ h^{\II^k}/ \ZZ D_k$.\end{remark}

We now construct the projective resolution of the object $\Delta_{\Box}\ZZ\in \Ab^{\Box}$. To this aim, consider homomorphisms 
$(d_k)_{\II^n}: \ZZ h^{\II^k}(\II^n)\to \ZZ h^{\II^{k-1}}(\II^n )$ assigning to each 
$f\in h^{\II^k}(\II^n)$ the sum 
$\sum\limits_{i=1}^k (-1)^i(f\circ\delta^{k,0}_i-f\circ\delta^{k,1}_i)$.
For an arbitrary $f\in D_k(\II^n)$, for $k\geq 1$, there are a morphism $g: \II^{k-1}\to \II^n$ and a number $j$ from range $1\leq j\leq k$ such that $f=g\circ\sigma^k_j$.
For an arbitrary $1\leq i\leq k$, we have the equality $f\delta^{k,\varepsilon}_i=g\sigma^k_j\delta^{k,\varepsilon}_i$.
By virtue of the formulas (\ref{rel3}) for $i<j$, we have
$g\sigma^k_j\delta^{k,\varepsilon}_i=
g\delta^{k-1,\varepsilon}_i\sigma^{k-1}_{j-1}\in D_{k-1}(\II^n)$.
Similarly, for $i>j$, $g\sigma^k_j\delta^{k,\varepsilon}_i\in D_{k-1}(\II^n)$.
If $i=j$, then $g\sigma^k_j\delta^{k,\varepsilon}_i=g$, and at the same time $f\circ\delta^{k,0}_i-f \circ\delta^{k,1}_i=g-g=0$. This implies that the homomorphisms $(d_k)_{\II^n}$ carry elements from $\ZZ D_k(\II^n)$ to elements from $\ZZ D_{k-1}(\II^n)$ .

Hence cocubical abelian groups $\ZZ h^{\II^k}/\ZZ D_k$
will form a chain complex whose differentials $\overline{d}_k$ have components
defined on cosets by subgroups
$\ZZ D_k(\II^n) \subseteq \ZZ h^{\II^k}(\II^n)$ by the formula
$$
 (\overline{d}_k)_{\II^n}(f+\ZZ D_k(\II^n))= (d_k)_{\II^n}(f)+\ZZ D_{k-1} (\II^n).
$$

Since $D_0(\II^n)=\emptyset$, then $\ZZ h^{\II^0}/\ZZ D_0=\ZZ h^{\II^0}$.

\begin{lemma}\label{exactseq}\cite[Lemma 5]{X2019}
For each $n\geq 0$ the complex of abelian groups
$$
 0 \leftarrow \ZZ h^{\II^0}/ \ZZ D_0(\II^n)
 \stackrel{(\overline{d}_1)_{\II^n}}\leftarrow
\ZZ h^{\II^1}/ \ZZ D_1(\II^n)
 \stackrel{(\overline{d}_2)_{\II^n}}\leftarrow
\ZZ h^{\II^2}/ \ZZ D_2(\II^n)
 \stackrel{(\overline{d}_3)_{\II^n}}\leftarrow \cdots.
$$
is isomorphic to the $C_*$ complex constructed in \cite{X2008},
consisting of Abelian groups and homomorphisms
$$
 0\leftarrow \ZZ\Box_+(\II^0,\II^n) \stackrel{d^+_1}\leftarrow
\ZZ\Box_+(\II^1,\II^n) \stackrel{d^+_2}\leftarrow \cdots
\stackrel{d^+_n}\leftarrow \ZZ\Box_+(\II^n,\II^n)\leftarrow 0.
$$
This means that its homology $H_k$ is equal to $0$ for $k>0$, and $H_0=\ZZ$.
\end{lemma}

We define a natural transformation
$\epsilon: \ZZ h^{\II^0}\to \Delta_{\Box}\ZZ$,
such that $\epsilon_{\II^n}: \ZZ\Box(\II^0,\II^n)\to \ZZ$ take on 
$x\in \Box(\II^0, \II^n)$ values $\epsilon_{\II^n}(x)=1$.

\begin{proposition}\label{normres}
Sequence of objects and natural transformations in $\Ab^{\Box}$
$$
0 \leftarrow \Delta_{\Box}\ZZ \stackrel{\epsilon}\leftarrow \ZZ h^{\II^0}/\ZZ D_0
\stackrel{d_1}\leftarrow \ZZ h^{\II^1}/\ZZ D_1
\stackrel{d_2}\leftarrow \ZZ h^{\II^2}/\ZZ D_2\leftarrow \cdots
$$
is the projective resolution of the diagram $\Delta_{\Box}\ZZ$.
\end{proposition}
{\sc Proof.} Lemma \ref{exactseq} implies that this sequence is exact.
 By the \ref{proj} lemma,
cubical Abelian groups $\ZZ h^{\II^k}/ \ZZ D_k$ will be
 projective objects of the category $\Ab^{\Box}$.
Therefore, this exact sequence will be a projective resolution.
\hfill$\Box$

\subsection{Homology of cubical objects in the AB4-category}

Let $\mA$ be an Abelian category.
Consider an arbitrary cubical object $F: \Box^{op}\to \mA$.
First, we introduce its unnormalized complex.
If $\mA=\Ab$ is the category of Abelian groups, then the normalized complex will consist of quotient groups $F(\II^k)$ with respect to subgroups generated by degenerate elements.

A quotient object in an Abelian category can be constructed as a subobject embedding cokernel. This construction will give a normalized complex in $\mA$.

Let $A$ and $B$ be objects of the Abelian category. For any morphism $f: A\to B$ we denote 
by $\Coker(f)$ its cokernel (respectively, by $\Ker(f)$ its kernel), and by 
$coker(f): B\to \Coker (f)$ the canonical projection onto the cokernel 
(resp. $ker(f): \Ker(f)\to A$ the canonical embedding of the kernel).

For any object $A$ of an abelian category and $n\geq 0$ denote by $A^n$ the coproduct of $n$ copies $A\oplus \cdots \oplus A$ of the object $A$. Let $in_i: A\to A^n$ be morphisms of the coproduct cone, $1\leq i\leq n$.
 For arbitrary two objects $A$, $B$ and morphisms $f_1, \ldots, f_n \in \mA(A,B)$ of the Abelian category, denote by $(f_1, \ldots, f_n): A^n\to B $ a morphism for which $(f_1, \ldots, f_n)\circ in_i= f_i$, for all $1\leq i\leq n$.
Let $\Coker(A^{n} \stackrel{(f_1, \ldots, f_n)}\longrightarrow B)$ be the cokernel of this morphism, and $coker(f_1, \ldots, f_n): B\to \Coker (A^{n} \stackrel{(f_1, \ldots, f_n)}\longrightarrow B)$
is the canonical projection. For example, in the case $\mA=\Ab$ the group $\Coker(f_1, \ldots, f_n)$ will be isomorphic to the factor group $B/(\Imm(f_1) + \ldots + \Imm(f_n))$.

Before defining the normalized complex corresponding to the cubical object $F$ of the abelian category, we note that the construction of the projective resolution for the diagram $\Delta_{\Box}\ZZ$, in Proposition \ref{normres}, was caused by the following idea:
The tensor product of this resolution by $F$ gives the complex
\begin{equation}\label{normcom}
0 \leftarrow (\ZZ h^{\II^0}/\ZZ D_0)\otimes F \stackrel{d_1}\leftarrow
(\ZZ h^{\II^1}/\ZZ D_1)\otimes F \stackrel{d_2}\leftarrow
(\ZZ h^{\II^2}/\ZZ D_2)\otimes F \leftarrow
\ldots,
\end{equation}
whose homology objects are isomorphic to $\coLim^{\Box^{op}}_n F$ 
by Proposition \ref{projcolim}.

{\it Normalized complex} of cubical object $F: \Box^{op}\to \mA$ consists of objects
$$
C^{\nd}_k(F)= \Coker (F(\II^{k-1})^k \stackrel{(F(\sigma^k_1), \ldots, F(\sigma^k_k))}
\longrightarrow F(\II^k)), \quad k\geq 0,
$$
given together with canonical morphisms 
$coker{(F(\sigma^k_1), \ldots, F(\sigma^k_k))}: F(\II^k)\to C^{\nd}_k(F)$.
We define the differentials of this complex as follows.
For any $i,j\in \{1, \ldots, k\}$ such that $i\neq j$, according to
relation (\ref{rel3}) there exists a pair of numbers $(i',j')$ for which the equality
$\sigma^k_j\delta^{k,\varepsilon}_i= \delta^{k-1,\varepsilon}_{i'}\sigma^{k-1}_{j'}$.
The pair $(i',j')$ is obtained by subtracting one from the larger number.
There is a commutative diagram
$$
\xymatrix{
F(\II^{k-1})\ar[d]_{F(\delta^{k-1,\varepsilon}_{i'})}
\ar[r]^{F(\sigma^k_j)} & F(\II^k) \ar[d]^{F(\delta^{k,\varepsilon}_i)}\\
F(\II^{k-2}) \ar[r]_{F(\sigma^{k-1}_{j'})} & F(\II^{k-1})
}
$$
Using this diagram and the equality
$F(\delta^{k,0}_i)F(\sigma^k_i)= F(\delta^{k,1}_i)F(\sigma^k_i)$,
construct morphisms 
$f_i: F(\II^{k-1})^k\to F(\II^{k-2})^{k-1}$,
making the right square of the following diagram commutative
$$
\xymatrix{
F(\II^{k-1}) \ar[d]|-{F(\delta^{k-1,0}_{i'})-F(\delta^{k-1,1}_{i'})}
\ar[r]^{in_j} & F(\II^{k-1})^k \ar@{-->}[d]^{f_i}
\ar[rrr]^{(F(\sigma^k_1),\ldots,F(\sigma^k_k))}
&&& F(\II^k) \ar[d]_{F(\delta^{k,0}_i)-F(\delta^{k,1}_i)}\\
F(\II^{k-2}) \ar[r]_{in_{j'}} & F(\II^{k-2})^{k-1}
\ar[rrr]_{~~(F(\sigma^{k-1}_1),\ldots,F(\sigma^{k-1}_{k-1}))}
&&& F(\II^{k-1})
}
$$
for each $i$ from the interval $1\leq i\leq k$.
To this purpose, given $i$, we define $f_i$ as the unique morphism of 
$F(\II^{k-1})^{\oplus k}\to F(\II^{k-2})^{\oplus(k-1)}$,
for which, for $j\not=i$ the composition $f_i\circ in_j$ equals
$in_{j'}\circ({F(\delta^{k-1,0}_{i'})-F(\delta^{k-1,1}_{i'})})$,
and for $j=i$ the composition of $f_i\circ in_j$ is equal to $0$.
The constructed morphism $f_i$ will satisfy the relations
$$
(F(\sigma^{k-1}_1),\ldots,F(\sigma^{k-1}_{k-1})) f_i in_j
=
(F(\delta^{k,0}_i)-F(\delta^{k,1}_i))(F(\sigma^{k}_1),
\ldots, F(\sigma^{k}_{k}))in_j.
$$
Since the morphism cone $in_j$ is separating, the right square is commutative.
Now it is easy to obtain a commutative diagram to which we can add a column with the morphism $d_k: C^{\nd}_k(F)\to C^{\nd}_{k-1}(F)$.
\begin{equation}\label{univ2}
\xymatrix{
  F(\II^{k-1})^k \ar[d]|-{\sum^k_{i=1}(-1)^i f_i}
\ar[rrr]^{(F(\sigma^k_1),\ldots,F(\sigma^k_k))}
&&& F(\II^k) \ar[rrrr]^{coker(F(\sigma^k_1),\ldots,F(\sigma^k_k))} \ar[d]|-{\sum^k_{i=1}(-1)^i(F(\delta^{k,0}_i)-F(\delta^{k,1}_i))}
& & & & C^{\nd}_k(F) \ar@{-->}[d]_{\exists!}^{d_k}\\
  F(\II^{k-2})^{k-1}
\ar[rrr]_{~~(F(\sigma^{k-1}_1),\ldots,F(\sigma^{k-1}_{k-1}))}
&&& F(\II^{k-1}) \ar[rrrr]_{coker(F(\sigma^{k-1}_1),\ldots,F(\sigma^{k-1}_{k-1}))} 
& & & & C^{\nd}_{k-1}(F)
}
\end{equation}
Since the compositions of row morphisms are zero and the left square of this diagram is commutative, there exists a unique morphism $d_k$ that makes the right square commutative.
Since the projections
 $coker(F(\sigma^k_1),\ldots,F(\sigma^k_k))$ are epimorphisms, and the middle column morphisms are complex, then $d_k\circ d_{k+1}=0$. Hence $(C^{\nd}_k(F), d_k)$  is a chain complex in the category  $\mA$.

\begin{theorem}\label{main1}
For an arbitrary cubical object $F$ in the AB4-category, the homology objects 
$H_k(C^{\nd}_*(F))$ are isomorphic to $\coLim^{\Box^{op}}_k F$, for all $k \geq 0$.
\end{theorem}
{\sc Proof.} Using the natural isomorphism (\ref{nattrans}) and the permutability 
of the functor $(-)\otimes F$ with colimits, we obtain an isomorphism of 
 the complexes $(C^{\nd}_k(F), d_k)$ and (\ref{normcom}).
By Proposition \ref{normres} the complex in $\Ab^{\Box}$ consisting of $\ZZ h^{\II^k}/\ZZ D_k$ is the projective resolution of the diagram $\Delta_{\Box}\ZZ$. Hence, using
 Proposition \ref{projcolim} we obtain $H_n(C^{\nd}_*(F))\cong \coLim^{\Box^{op}}_n F$.\hfill$\Box$

\begin{example}
Consider $\mA= \Ab$. Let $F: \Box^{op}\to \Ab$ be a cubical Abelian group. In this case, 
the complex whose homology groups are isomorphic to $\coLim^{\Box^{op}}_n F$ 
will consist of the factor groups $C^N_k(F)= F(\II^k)/ \sum^ k_{i=1} \Imm F(\sigma^k_i)$.
Therefore, the differential $d^N_k$ associates with the coset containing $a\in F(\II^k)$ 
the coset of the element $d_k(a)\in F(\II^{k-1})$.
\end{example}

\subsection{Cohomology of cocubical objects}

Let $A$ be an object of the Abelian category $\mA$. Denote by $pr_i: A^n\to A$ the morphisms of the product cone. For arbitrary objects $A, B\in \mA$ and morphisms $f_1, \ldots, f_n \in \mA(B, A)$ there is a morphism which we denote by $(f_1, \ldots, f_n)^*: B\to A^n$ satisfying the relations $f_i= pr_i\circ(f_1, \ldots, f_n)^*$ for all $1\leq i\leq n$.

The functor $F: \Box \to \mA$ is called a cocubical object in the category $\mA$.
It can be associated with the cubical object $F^{op}: \Box^{op}\to \mA^{op}$ in the category $\mA^{op}$. If $\mA$ is an AB4-category, then the homology for $F^{op}$ is defined. The objects dual to the homology of the cubical objects $F^{op}$ in $\mA^{op}$ are called the cohomology of the cocubical object $F: \Box\to \mA$.
Using the principle of duality, we obtain a normalized complex for cohomology. It consists of objects

$$
C^k_{\nd}(F)= Ker ( F(\II^k) \stackrel{(F(\sigma^k_1), \ldots, F(\sigma^k_k))^*}
\longrightarrow F(\II^{k-1})^k).
$$
The diagram (\ref{univ2}) in the category $\mA^{op}$ contains a commutative square

\begin{equation}\label{subsquare}
\xymatrix{
F(\II^k) &&&& C^k_N(F) \ar[llll]_{ker(F(\sigma^k_1), \ldots, F(\sigma^k_k))^*}\\
F(\II^{k-1}) \ar[u]|-{\sum^k_{i=1}(-1)^i(F(\delta^{k,0}_i)-F(\delta^{k,1}_i))}
&&&& C^{k-1}_N(F) \ar[u]_{d^{k-1}_N}
\ar[llll]^{ker(F(\sigma^{k-1}_1), \ldots, F(\sigma^{k-1}_{k-1}))^*}
}
\end{equation}
leading to the definition of differentials $d^k_N: C^k_N(F)\to C^{k+1}_N(F)$ for $k\geq 0$. (For $k<0$ the differentials and objects are equal to $0$.)

The following statement is obtained from Theorem \ref{main1},
if we replace the category $\mA$ in it with the dual category.

\begin{corollary}\label{maincor1}
For an arbitrary cocubical object $F$ in the AB4*-category, the cohomology objects $H^k(C^*_{\nd}(F))$ are isomorphic to $\Lim^k_{\Box} F$, for all $k\geq 0$.
\end{corollary}

\begin{example}
Consider $\mA=\Ab$. For a cocubical abelian group $F: \Box\to \Ab$ the cochain normalized complex will consist of the abelian groups
$$
C^k_N(F)= \bigcap\limits^k_{i=1}Ker(F(\sigma^k_i))= \{a \in F(\II^k)~|~ 
(\forall i\in \{1, \ldots, k)\}) F(\sigma^k_i)(a)=0\}
$$ 
and homomorphisms $d^k(a)= 
{\sum^{k+1}_{i=1}(-1)^i(F(\delta^{k+1,0}_i)-F(\delta^{k+1,1}_i))}$ 
defined as differentials of the complex corresponding to the cocubical abelian group $F$. 
Since the diagram (\ref{subsquare}) is commutative, these differentials carry elements from $C^k_N(F)$ to elements from $C^{k+1}_N(F)$.
\end{example}

\section{Homology of cubical sets with coefficients in systems}

Let $\mA$ be an AB4-category. 

The category of cubes can be considered as a full subcategory of the category 
$\Set^{\Box^{op}}$ as an image of the Yoneda embedding 
$h^{\Box}: \Box \xrightarrow{\subseteq} \Set^{\Box^{op}}$.

Consider a cubical set $X\in \Set^{\Box^{op}}$. Left fiber of $h^{\Box}$ over $X$  denote by $\Box/X$.
Let $Q_X: \Box/X\to \Box$ be a forgetting functor of the left fiber.

A {\it contravariant system of objects} in $\mA$ on $X$ is a diagram 
$F: (\Box/X)^{op}\to \mA$. 
The homology objects of $C^{\nd}_*(Lan^{Q^{op}_X}F)$ are called {\it the homology objects
 of the cubical set $X$ with coefficients in $F$}.
Here we construct a chain complex whose homology objects are isomorphic to the homology objects of the cubical set $X$ with coefficients in $F$.

\subsection{Construction of a normalized complex}\label{buildcom}

The set of objects in the category $(\Box/X)^{op}$ can be considered as a disjoint union $\coprod\limits_{n\geq 0}X_n$.
Morphisms in it, from $(x,m)$ to $(y,n)$, can be considered as triples 
$x\stackrel{\alpha}\to y$ such that $\alpha\in \Box(\II^n, \II^m)$ and 
$X(\alpha)(x)=y$.
Consider a cubical object 
$(C_n(X,F), d^{n,\varepsilon}_i, s^n_i):= Lan^{Q^{op}_X}F$ in the category 
$\mA$. According to the proposition \ref{contralan}(\ref{conten1}) its objects are equal to $C_n(X,F)=\bigoplus\limits_{x\in X_n}F(x)$.
It follows from Proposition \ref{contralan}(\ref{conten1}) that the boundary 
operators $d^{n,\varepsilon}_i: C_{n}(X,F)\to C_{n-1}(X, F)$ are determined by the commutativity condition for diagrams

\begin{equation*}
\xymatrix{
\bigoplus\limits_{x\in X_n}F(x) \ar[rrr]^{d_i^{n,\varepsilon}}
 &&& \bigoplus\limits_{x\in X_{n-1}}F(x)\\
F(x) \ar[u]^(.4){in_x} 
\ar[rrr]^(.45){F(\delta_i^{n,\varepsilon}:x\to X(\delta^{n,\varepsilon}_i)x)} 
&&&
F(X(\delta_i^{n,\varepsilon})x)\ar[u]_(.35){in_{X(\delta^{n,\varepsilon}_i)x}}
}
\end{equation*}

It also follows from the proposition \ref{contralan}(\ref{conten1}) that for the degeneration operators $s^n_i: C_{n-1}(X,F)\to C_{n}(X,F)$ are commutative diagrams

\begin{equation}\label{comdeg}
\xymatrix{
\bigoplus\limits_{x\in X_{n-1}}F(x) \ar[rrr]^{s_i^{n}} &&&
	\bigoplus\limits_{x\in X_{n}}F(x)\\
F(x)\ar[u]^(.4){in_x}
\ar[rrr]^(.45){F(\sigma_i^{n}:x\to X(\sigma^{n}_i)x)} &&&
F(X(\sigma^{n}_i)x) \ar[u]_(.35){in_{X(\sigma^n_i)x}}
}
\end{equation}

Let $d_n=\sum\limits^n_{i=1}(-1)^i(d^{n,0}_i-d^{n,1}_i)$.
Consider a sequence of objects

$$
C^{\nd}_n(X,F)= \Coker\left((\oplus_{x\in X_{n-1}}F(x))^{\oplus n}
\xrightarrow{(s^n_1, \ldots, s^n_n)}
\oplus_{x\in X_n}F(x)\right)
$$
together with canonical projections
$\bigoplus\limits_{x\in X_n}F(x)\stackrel{pr_n}\longrightarrow C^{\nd}_n(X,F)$.

From the commutativity of the diagram (\ref{univ2}) for an arbitrary cubical object in the Abelian category, it follows that the composition $pr_{n-1}\circ d_n \circ (s^n_1, \ldots, s^n_n)$ is equal to zero.
Since $pr_n\circ (s^n_1, \ldots, s^n_n)=0$, and $pr_n$ as a projection onto the cokernel has the universality property, there exists a unique morphism $d^N_n: C^{\nd}_n( X,F)\to C^{\nd}_{n-1}(X,F)$ satisfying $d^N_n\circ pr_n= pr_{n-1}\circ d_n$.
Thus, the normalized complex $(C^{\nd}_n(X,F), d^{\nd}_n)_{n\geq 0}$ is constructed, whose homology objects are taken as the definition of the homology objects for the cubical set $X$ with coefficients in the contravariant system
 $F: (\Box/X)^{op} \to \mA$.

\begin{example}
Let $X$ be a cubical set and let $G: \Box/X)^{op}\to \Ab$ be a contravariant system
 of abelian groups. The complex obtained from the cubical Abelian group $Lan^{Q^{op}_X}G: \Box^{op}\to \Ab$ has boundary operators acting as
$$
d^{k,\varepsilon}_i(x\in X_k, g\in G(x))= (X(\delta^{k,\varepsilon}_i)x, 
G(\delta^{k,\varepsilon}_i: x\to X(\delta^{k,\varepsilon}_i)x)).
$$ 
It consists of abelian groups
$C_k(X, G)= \oplus_{x\in X_k}G(x)$ and differentials 
$d_k= \sum^k_{i=1}(-1)^i(d^{k,0}_i- d^{k,1}_i)$. 
 For each $x\in X_{k-1}$, the degeneracy operator  $s^k_i= Lan^{Q^{op}_X}G: C_{k-1}(X, G)\to C_k(X,G)$ in this cubical abelian group acts as
  $(x, g\in G(x)) \mapsto (X(\sigma^k_i)x, G(x\stackrel{\sigma^k_i}\to X(\sigma^k_i)x)(g))$.
The normalized complex will consist of quotient groups
$$
C^N_k(X,G)=
\oplus_{x\in X_k}G(x) /{\sum^k_{i=1}\Imm s^k_i}. 
$$
Differentials $d^N_k$ are defined by differentials $d_k$. 
\end{example}

\begin{example}
Consider a covariant system $G: \Box/X\to \Ab$ in $\Ab$. 
The objects of the category $\Box/X$ 
are elements $x\in \coprod\limits_{n\in \NN}X_n$, and morphisms in $\Box/X$ are 
defined as triples $y\xrightarrow{\alpha} x$ for which $x \xrightarrow{\alpha} y$ are morphisms in $(\Box/X)^{op}$, which means $X(\alpha)(x)=y$.
 
The cocubical Abelian group $Ran_{Q_X}G: \Box \to \Ab$ gives a cochain complex consisting of Abelian groups $C^k(X,G)= \prod_{x\in X_k}G(x)$.
Elements of $C^k(X,G)$ are given as functions $\varphi: X_k\to \bigcup_{x\in X_k}G(x)$ such that $\varphi(x)\in G(x )$ for all $x\in X_k$.
Its differentials are $d^k= \sum^k_{i=1}(-1)^i(d^i_{k,0}- d^i_{k,1})$, where $d^i_{ k,\varepsilon}: C^{k-1}(X,G)\to C^{k}(X,G)$ are coboundary operators defined by the following formula
for $\varphi\in C^{k-1}(X,G)$:
$$
d^i_{k,\varepsilon}(\varphi)(x)= G(\delta^{k, \varepsilon}_i: X(\delta^{k, \varepsilon}_i)x \to x)\varphi(X (\delta^{k, \varepsilon}_i)x).
$$

The cocubical Abelian group $Ran_{Q_X}G$ also gives the degeneration operations $s^i_k: C^k(X, G)\to C^{k-1}(X,G)$ of this complex defined by the formula
$$
s^i_k(\varphi)(x)= G(X(\sigma^k_i)x\xrightarrow{\sigma^k_i} x)\varphi(X(\sigma^k_i)x)
$$
So, the normalized complex consists of the groups $C^k_N(X, G)= \bigcap^k_{i=1} Ker(s^i_k)$, and its differentials are defined as $d^k_N(\varphi)= d^k( \varphi)$.
\end{example}

\subsection{Homology of a cubical set and satellites of the colimit functor}

Consider the normalized complex $C^{\nd}_*(X,F)=(C^{\nd}_n(X,F), d^{\nd}_n)$
 built in Sub\-sec\-ti\-on \ref{buildcom} for the functor $F: \Box/X)^{op}\to \mA$, 
 into the AB4-category $\mA$. 
We have established that the homology $H_n(C^{\nd}_*(X,F))$ of this complex, $n\geq 0$, is equal to the homology of the cubical set $X$ with coefficients in $F$.
 Denote the homology of $X$ with coefficients in $F$ by $H_n(X, F)$, $n\geq 0$.    

\begin{theorem}\label{main2}
Let $X$ be a cubical set, and let $F$ be a contravariant system of objects in the AB4-category 
$\mA$ on $X$. Then $H_n(X,F) \cong \coLim^{(\Box/X)^{op}}_n F$, 
for all $n\geq 0$.
\end{theorem}
{\sc Proof.} We first constructed the cubical object in the Abelian category as the left Kan extension $Lan^{Q_X^{op}}F: \Box^{op}\to \mA$.
The complex $C^{\nd}_*(X,F)$ was constructed as a normalized complex of the cubical object $Lan^{Q_X^{op}}F$.
Since for each $\II^k\in \Box$ each connected component of the category $\II^k/Q_{X}$ has 
an initial object, the values of $\coLim^{(\Box/X)^{op}}_n F$ are isomorphic to 
$\coLim^{\Box^{op}}_n Lan^{Q_X^{op}}F$, which is isomorphic to 
$H_n(C^{\nd}_*( Lan^{Q_X^{op}}F))$.
\hfill$\Box$

\begin{corollary}\label{main2coh}
Let $X$ be a cubical set, $F$ be a covariant system of objects in the AB4-category $\mA$ on $X$. Then there are isomorphisms $H^n(C_{\nd}^*(X,F)) \cong \Lim_{\Box/X}^n F$, for all $n\geq 0$.
\end{corollary}

\subsection{Homology invariance upon passing to the direct image}

For every cubical set morphism $f: X\to Y$, the functor $\Box/f: \Box/X\to \Box/Y$ is defined. It assigns to each $\widetilde{x}: \Box^n\to X$ object of $\Box/X$ an object of $\Box/Y$ that is equal to the composition $\Box^n\stackrel{\widetilde{x} }\to X\stackrel{f}\to Y$. Morphisms 
in $\Box/X$ are commutative triangles with sides $(\Box^{p_1}\stackrel{\Box(\alpha)}\to \Box^{p_2}, \Box^{p_1}\stackrel{\widetilde{x_1}}\to X,
\Box^{p_2}\stackrel{\widetilde{x_2}}\to X)$ and they go to
  $(\Box^{p_1}\stackrel{\Box(\alpha)}\to \Box^{p_2},
   \Box^{p_1}\stackrel{f\widetilde{x_1}}\to X,
\Box^{p_2}\stackrel{f\widetilde{x_2}}\to Y)$.

By Theorem \ref{lanhomol}, 
$\coLim^{(\Box/X)^{op}}_n F \stackrel{\cong}\to \coLim^{(\Box/Y)^{op}}_n f_* F$.
Applying Theorem \ref{main2}, we get

\begin{corollary}\label{dirhomol}
Let $f: X\to Y$ be a morphism of cubical sets.
Then for any functor $F: (\Box/X)^{op}\to \mA$
into the AB4-category, there are natural isomorphisms
$H_n(X, F)\stackrel{\cong}\to H_n(Y, f_* F)$, for all $n\geq 0$.
\end{corollary}

\subsection{Inverse image homology preservation criterion}

Let $f: X\to Y$ be a morphism of cubical sets.
For each cube $\tilde{y}: \Box^n\to Y$, its inverse fiber 
$\overleftarrow{f}({y})$ is defined by pullback (\ref{invimage}). The proposition \ref{critiso} and the theorem \ref{main2} imply the following assertion:

\begin{corollary}\label{critcube}
The following properties of the cubical set morphism $f: X\to Y$ are equivalent.
\begin{enumerate}
\item For all $y\in Y$ and $n\geq 0$ the groups $H_n(\overleftarrow{f}(y), \Delta\ZZ)$ are isomorphic to the homology groups of the cubical point $H_n(\Box^0, \Delta\ZZ)=\begin{cases}
\ZZ, & \text{ if $n=0$;}\\
0, & \text{ if $n>0$.}
\end{cases}$
\item The canonical homomorphisms of Abelian groups $H_n(X, f^*F) \to H_n(Y,F)$ are isomorphisms for every functor $F: (\Box/Y)^{op}\to \Ab$, for all $n\geq 0$.
\item The canonical morphisms $H_n(X, f^*F) \to H_n(Y,F)$ are isomorphisms for any functor $F: (\Box/Y)^{op}\to \mA$ into an arbitrary AB4-category $\mA$, 
for all $n\geq 0$.
\end{enumerate}
\end{corollary}

\subsection{Spectral sequence for the colimit homologies of cubical sets}

From the theorem \ref{main2}, substituting the category $\Box$ instead of the small category $\mD$ into the sentence \ref{covering}, we arrive at a spectral sequence 
converging to the homology of the colimit of cubical sets.
A similar result for semi-cubical sets and $\mA= \Ab$ was obtained in \cite{X2008}.

\begin{corollary}\label{coveringcor}
Let $J$ be a small category, and let $\{X^i\}_{i\in J}$ be a cubical set diagram such that
\begin{equation}\label{c3cor}
\coLim_q^J\{\ZZ(X^i_n)\}_{i\in J}= 0, \mbox{ for any } n\geq 0 \mbox{ and } q>0.
\end{equation}
Let $\lambda_i: X^i\rightarrow \coLim^J \{X^i\}_{i\in J}$ be the cone of colimit morphisms of cubical sets.
Then for any AB4-category $\mA$ and any functor 
$F: (\Box/\coLim^J \{X^i\} )^{op}\to \mA$ 
there exists a spectral sequence of the first quadrant
$$
E^2_{p,q}= \coLim_p^{J}\{H_q(X^i, \lambda^*_i F)\}_{i\in J} \Rightarrow H_{p+q}(\coLim ^J\{X^i\}, F).
$$
\end{corollary}

The condition (\ref{c3cor}) is satisfied, for example, for a locally directed cover 
of a cubical set. In this case, we get a generalization of \cite[Corollary 5.3]{X2008}, proved there for semi-cubical sets.

\subsection{Cubical homology of small categories}

We introduce the cubical homology of a small category as the cubical homology of its cubical nerve and prove that these homology are isomorphic to the small category homology defined in \S2.

  A cubical nerve of the small category $\mC$ \cite{jar2006} is a cubical set 
  $\nr^{\Box}\mC:= \Cat(-, \mC)|_{\Box}: \Box^{op}\to \Set$. It is defined as the restriction of the functor $\Cat(-, \mC)$ to the category of cubes. The set of its $n$-dimensional cubes is $\nr^{\Box}_n(\mC)= \Cat(\II^n,\mC)$ for all $n\geq 0$.
  
 The category $\Box/\nr^{\Box}\mC$ is isomorphic to the comma category $\Box/\mC$ whose objects are the morphisms $\II^n\to \mC$, $n\geq 0$, and whose morphisms from $\II^m\xrightarrow{f} \mC $ to $\II^n\xrightarrow{g} \mC $ are commutative triangles
\begin{equation}\label{morphoverc}
\xymatrix{
\II^m\ar[rd]_{\theta} \ar[rr]^f && \mC\\
& \II^n \ar[ru]_g
}
\end{equation}
where $\theta$ is a morphism of the category $\Box$.
We will identify the categories $\Box/\nr^{\Box}\mC$ and $\Box/\mC$.

We define the functor $\partial: \Box/\mC\to \mC$ as assigning to each box 
$x: \II^n\to \mC$ an object $x(1, \cdots, 1)$.
Recall that a functor $S:\mC\to \mD$ is said to be aspherical if its left fibers $S/d$
 are contractible \cite{jar2006}. In this case $\coLim^{\mC^{op}}_n FS \cong \coLim^{\mD^{op}}_n F$.
The functor $\partial$ is aspherical \cite[Lemma 23]{X2019}.

Let $H_n(\nr^{\Box}\mC, F\circ\partial^{op})$ be the $n$th cubical homology 
of a cubical nerve. Theorem \ref{main2} and the asphericity of the functor $\partial$ imply
  
\begin{corollary}\label{homolcatcub}
Let $\mC$ be a small category and $\mA$ an AB4-category. Then for every functor $F: \mC^{op}\to \mA$ the cubical homology objects $H_n(\nr^{\Box}\mC, F\circ\partial^{op})$ are isomorphic to $\coLim_n^{ \mC^{op}}F$, for all $n\geq 0$.
\end{corollary}

This assertion was proved in \cite[Corollary 24]{X2019} for the case $\mA= \Ab$.

\subsection{Cubical Baues-Wirsching cohomology}

Let $\mA$ be an AB4*-category.
Denote by $H^n_{BW}(\mC,F)$, $n\geq 0$, the $n$th Baues-Wirsching cohomology of the small category with coefficients in the natural system of objects in $\mA$ \cite{bau1985}, \cite{gal2012}.

We introduce the functor $\frak{d}: \Box/\mC\to \fF\mC$.
For this purpose, denote by $0^n\in \II^n$, $n>0$, the point of the cube, all of whose coordinates are equal to $0$, and by $1^n$, the point of the cube, all of whose coordinates are equal to $1$. If $n=0$, then these points coincide, i.e. $0^0=1^0$.

We define the functor $\frak{d}: \Box/\mC\to \fF\mC$ on the objects $\II^n \xrightarrow{g} \mC$ as $\frak{d}(g)= (g( 0^n)\xrightarrow{g(0^n\leq 1^n)} g(1^n))$, $n\geq 0$.
For each morphism $\theta: f\to g$ from $\Box/\mC$ shown in the diagram (\ref{morphoverc}), the equality $g\theta= f$ holds.
We assign to it a morphism $\frak{d}(f)\to \frak{d}(g)$ of the category $\fF\mC$ for which the following diagram is commutative:

$$
\xymatrix{
f(1^m)=g\theta(1^m) \ar[rr]^(.6){g(\theta(1^m)\leq 1^n)}  && g(1^n)\\
f(0^m)=g\theta(0^m)
\ar[u]^{f(0^m\leq 1^m)}   && g(1^n)\ar[ll]^(.4){g(0^n\leq \theta(0^m))}
\ar[u]_{g(0^n\leq 1^n)}
}
$$ 

\begin{corollary}\label{homolbwcub}
Let $\mC$ be a small category, $G: \fF\mC\to \mA$ be a natural system on $\mC$.
Then $H^n_{BW}(\mC, G)\cong H^n(\nr^{\Box}\mC, G\circ\frak{d})$,
for all $n\geq 0$.
\end{corollary}
{\sc Proof.} Consider an arbitrary morphism $\alpha\in Mor\mC= \Ob\fF\mC$. Let us 
prove that the left fiber $\frak{d}/\alpha$ is contractible. Let $Id\langle\alpha\rangle$ 
be a category whose objects consist of decompositions of the morphism $\alpha$ represented as a pair of morphisms $\dom\alpha \xrightarrow{\beta} c \xrightarrow{\gamma}\cod\alpha$ category $\mC$ such that $\gamma\beta=
\alpha$.
Morphisms between the decompositions of $\dom\alpha \xrightarrow{\beta} c \xrightarrow{\gamma}\cod\alpha$ and $\dom\alpha \xrightarrow{\beta'} c' \xrightarrow{\gamma'}\cod \alpha$ are defined using the morphisms $c\xrightarrow{\xi} c'$ making the following diagram commutative

$$
\xymatrix{
\dom\alpha\ar[rd]_{\beta'} \ar[r]^{\beta} & c\ar[d]_{\xi} \ar[r]^{\gamma} & \cod{\alpha}\\
& c' \ar[ru]_{\gamma'}\\
}
$$

The category $Id\langle\alpha\rangle$ has an initial object $\dom\alpha \xrightarrow{1_{\dom\alpha}} \dom\alpha \xrightarrow{\alpha}\cod\alpha$ and a final object $\dom \alpha \xrightarrow{\alpha} \cod\alpha
\xrightarrow{1_{\cod\alpha}}\cod\alpha$, and hence its nerve is contractible.

Each object $\frak{d}(\II^n\xrightarrow{g} \mC) \xrightarrow{(u_0, u_1)} \alpha$
 of the category $\frak{d}/\alpha$ is defined by a pair consisting of the functor $\II^n\xrightarrow{g}\mC$ and the morphism $(u_0, u_1)$  of the category $\fF\mC$
 such that the following diagram $\in \mC$ is commutative:

$$
\xymatrix{
g(1^n)\ar[rr]^{u_1} && \cod\alpha\\
g(0^n)\ar[u]^{g(0^n\leq 1^n)} && \dom\alpha \ar[ll]_{u_0} \ar[u]^{\alpha}
}
$$
Associate with each such pair a functor $\II^n\to Id_{\mC}\langle\alpha\rangle$ taking each point $x\in \II^n$ to an object of the category $Id_{\mC}\langle\alpha\rangle$ 
equal to decomposition
$$
 \dom\alpha \xrightarrow{\beta_x} g(x) \xrightarrow{\gamma_x} \cod\alpha  
$$
where $\beta_x$ is equal to composition of morphisms $\dom\alpha \to g(0^n)\to g(x)$, 
and $\gamma_x$ is equal to composition $g(x)\to g(1^n)\to \cod\alpha$.

To each morphism $x\leq x'$ of the category $\II^n$ we associate a morphism in $Id_{\mC}\langle\alpha\rangle$ given by the diagram

$$
\xymatrix{
\dom\alpha\ar[rd]_{\beta_{x'}} \ar[r]^{\beta_x} & g(x)\ar[d]|-{g(x\leq x')} \ar[r]^{\gamma_x} & \cod{\alpha}\\
& g(x') \ar[ru]_{\gamma_{x'}}\\
}
$$
This leads to a bijection between $\Ob(\frak{d}/\alpha)$ and $\Ob(\Box/Id_{\mC}\langle\alpha\rangle)$.
Morphisms between objects from $\frak{d}/\alpha$ are defined using morphisms of the category $\Box$. The same is true for morphisms of the category $\Box/Id_{\mC}\langle\alpha\rangle$.
Hence, each morphism of the category $\frak{d}/\alpha$ must be associated with a morphism of the category $\Box/Id_{\mC}\langle\alpha\rangle$ having a morphism of cubes equal to the morphism of cubes between the corresponding objects.

We obtain an isomorphism of the categories $\frak{d}/\alpha\cong \Box/Id_{\mC}\langle\alpha\rangle$.
Moreover, the simplicial nerve of the category $Id_{\mC}\langle\alpha\rangle$ is contractible, and hence $\Box/Id_{\mC}\langle\alpha\rangle$ has a contractible nerve.
Consequently, the simplicial nerve of the category $\frak{d}/\alpha$ is contractible, whence for every functor $G: \fF\mC\to \Ab$ there is a natural isomorphism $H^n_{BW}(\mC, G)\cong H^n(\nr^{\Box}\mC, F\frak{d})$.
\hfill$\Box$

\section{Local homology and cohomology of cubical sets}

Consider the homology of cubical sets with coefficients in local systems of objects in the AB4-category.
We will prove that they do not depend on degenerate cubes and objects of the local system defined on degenerate cubes.
Then we consider the cohomology of cubical sets with coefficients in local systems and prove that the weak equivalence of cubical sets gives an isomorphism of cohomology with coefficients in local systems of abelian groups.

\subsection{Homology of cubical sets with coefficients in local systems}

A contravariant system on a cubical set is called local if its values on morphisms are isomorphisms. In this subsection, we first notice that the homologies of the standard cube with coefficients in the local system are isomorphic to the homologies of a point. 
Then we prove that the homology of a cubical set with coefficients in the local system is isomorphic to the homology of the complex with $n$th chain objects $\oplus_{x\in X^N_n} F(x)$, for all $n\geq 0$.

\begin{corollary}\label{homolcube}
Let $F$ be a local system on the cube $\Box^n$. Then 
$$
H_q(\Box^n, F)=\begin{cases}
F(1_{\II^n}), & \text{ if $q=0$;}\\
0, & \text{ if $q>0$.}
\end{cases}
$$
\end{corollary}
{\sc Proof.}  
To an arbitrary small category $\mC$ and a diagram $G: \mC\to \mA$ consisting of isomorphisms, there corresponds a diagram $G^{-1}$ obtained from $G$ by inverting its morphisms.
According to \cite[Appendix II, Proposition 4.4]{gab1967}, there are isomorphisms $\coLim^{\mC}_k G \cong \coLim^{\mC^{op}}_k G^{-1}$, for all $k\geq 0$.
The category $\Box/\Box^n$ has a final object $\tilde{1}_{\II^n}: \Box^n \to \Box^n$.
Hence $H_k(\Box^n, F) \cong \coLim^{(\Box/\Box^n)^{op}}_k F \cong \coLim^{\Box/\Box^n}_k F ^{-1}$, whence the proof follows.
\hfill$\Box$

\smallskip

Let us construct a (reduced) complex for finding homology of cubical sets
 with coefficients in local systems.

Let $X$ be a cubical set. Denote by $X^N_n$ the set of non-degenerate $n$-dimensional cubes, for $n\geq 0$.

\begin{theorem}\label{comloc}
Let $\mA$ be an AB4-category, $X$ a cubical set, $F: (\Box/X)^{op}\to \mA$ a local system on $X$.
  Then the objects of its normalized complex $C^{\nd}_n(X,F)$ are isomorphic to 
  $\bigoplus\limits_{x\in X^{\nd}_n}F(x)$, and the differentials $d^{\nd}_n: \bigoplus\limits_{x\in X^{\nd}_n}F(x)\to  \bigoplus\limits_{x\in X^N_{n-1}}F(x)$ consist of morphisms making the following diagram commutative
  
\begin{equation}\label{morn}
\xymatrix{
\bigoplus\limits_{x\in X_n}F(x)\ar[d]_{pr_n} \ar[rr]^{d_n} &&
\bigoplus\limits_{x\in X_{n-1}}F(x)\ar[d]^{pr_{n-1}}\\
\bigoplus\limits_{x\in X^N_n}F(x) \ar[rr]_{d^N_n}
&& \bigoplus\limits_{x\in X^N_{n-1}}F(x)
}
\end{equation}
Here $d_n=\sum\limits_{i=1}^{n}(-1)^i(d^{n,0}_i-d^{n,1}_i)$ 
are the differentials of the complex corresponding to the cubical object
$Lan^{Q^{op}_X}F$.
\end{theorem}
{\sc Proof.}
We have defined chain objects of the normalized complex by the formula
$$
C^{\nd}_n(X,F)=\Coker\left((\oplus_{x\in X_{n-1}}F(x))^{\oplus n}
\xrightarrow{(s^n_1, \ldots, s^n_n)} \oplus_{x\in X_n}F(x)\right)
$$

It suffices to prove that for every $n\geq 0$, the object $\oplus_{x\in X^N_n}F(x)$ together with the morphism $\oplus_{x\in X_n}F(x)\stackrel{pr_n }\to \oplus_{x\in X^N_n}F(x)$, defined as a projection onto the direct summand, is the cokernel of the morphism $(s^n_1, \ldots, s^n_n)$.
Morphisms $d^N_n$ will be defined using the universality property of the cokernel functor.

The morphisms $s^n_i: \oplus_{x\in X_{n-1}} F(x)\to \oplus_{x\in X_{n}} F(x)$ are defined using the commutativity property of the diagram (\ref {comdeg}) as satisfying the relations

$$
s^n_i\circ in_x = in_{X(\sigma^n_i)(x)}\circ F(x\xrightarrow{\sigma^n_i} X(\sigma^n_i)(x))
$$
for all $n\geq 1, 1\leq i\leq n$, and $x\in X_{n-1}$.

The image of the degeneration operator $s^n_i$ is equal to the sum of the images of the morphisms $s^n_i in_x$ over all $x\in X_{n-1}$.
Hence, since the morphisms $F(x\xrightarrow{\sigma^n_i} X(\sigma^n_i)(x))$ are isomorphisms, the images of the operators $s^n_i$ are equal to the direct sum of the objects $F(X(\sigma^n_i)x)$. This direct sum is equal to $\oplus_{x\in D_n X}F(x)$,
where $D_n X\subseteq X_n$ is a subset of degenerate $n$-dimensional cubes, which implies that the cokernel of the morphism $(s^n_1, \ldots, s^n_n)$ is equal to $\oplus_{x\in X^ N_n}F(x)$.

This could completes the proof.
But we have applied some intuitive terms without substantiating them.
We need to prove that the cokernel of the morphism $(s^n_1, \ldots, s^n_n)$ is equal to $\oplus_{x\in X^N_n} F(x)$.

For this purpose, consider an arbitrary object $A\in \mA$ and a morphism $\gamma: \oplus_{x\in X_n}F(x)\to A$.
For them there are implications

$$
\aligned
\gamma\circ(s^n_1, \ldots, s^n_n)=0
&\Leftrightarrow  (\forall i\in\{1,\ldots, n\})
\gamma\circ s^n_i=0 \\
&\Leftrightarrow (\forall i\in\{1,\ldots, n\})
(\forall x\in X_{n-1}) \gamma\circ s^n_i\circ in_x =0,
\endaligned
$$
where $in_x: F(x) \to \oplus_{x\in X_{n-1}}F(x)$ are 
the canonical coproduct morphisms.
Using the commutative diagram (\ref{comdeg}), we get
$$
\gamma\circ(s^n_1, \ldots, s^n_n)=0
 \Leftrightarrow \gamma\circ in_{X(\sigma^n_i)x}\circ
 F(\sigma^n_i: x\to X(\sigma^n_i)x)=0.
$$

After that, we use the fact that $F$ is a local system.
This gives the invertibility of the morphisms $F(\sigma^n_i: x\to X(\sigma^n_i)x)$ and, at the same time, the logical equivalences
$$
\aligned
\gamma\circ(s^n_1, \ldots, & s^n_n)=0 \\
 &\Leftrightarrow (\forall i\in\{1, \ldots, n\})(\forall x\in X_{n-1}) \gamma\circ in_{X(\sigma^n_i)x}=0 \\
 &\Leftrightarrow (\forall x\in X_n\setminus X^N_n) \gamma\circ in_x =0.
\endaligned
$$

Since $(\forall x\in X_n\setminus X^N_n) pr_n\circ in_x=0$, from here we get $pr_n\circ(s^n_1, \ldots, s^n_n)=0$.

Moreover, it is easy to see that for any $\gamma: \oplus_{x\in X_n}F(x)\to A$ satisfying the condition $\gamma\circ (s^n_1, \ldots, s^n_n)= 0$, there exists a morphism $\widetilde{\gamma}:\oplus_{x\in X^N_n}F(x)\to A$ such that $\widetilde{\gamma}\circ pr_n=\gamma$.
For this purpose, we can take $\widetilde{\gamma}= \gamma\circ in_n$, where $in_n: \oplus_{x\in X^N_n}F(x) \to \oplus_{x\in X_n}F( x)$ - canonical direct summand injection.
To prove the equality $\widetilde{\gamma}\circ pr_n=\gamma$, it suffices to consider the morphisms $in_x$, for all $x\in X_n$, and consider the compositions $\gamma\circ in_n \circ pr_n \circ in_x$ in two cases - when $x\in X_n\setminus X^N_n$ and when $x\in X^N_n$.
In all cases, we get $\gamma\circ in_n\circ pr_n\circ in_x= \gamma\circ \in_x$.
By virtue of the separating property of the coproduct morphism cone, we obtain $\gamma\circ in_n=\gamma$.

Since $pr_n$ is an epimorphism, $\widetilde{\gamma}$ will be unique.
Hence the pair $(\oplus_{x\in X^N_n}F(x), pr_n)$ is the cokernel of the morphism $(s^n_1, \ldots, s^n_n)$, and hence it is isomorphic to $C^{\nd }_n(X,F)$.

The morphism $d^N_n$ was constructed above for an arbitrary functor $F$ as a morphism between cokernels. In the case when $F$ is a local system, it will be a morphism making the diagram (\ref{morn}) commutative.
Since $pr_n$ is an epimorphism, $d^N_n$ will be unique.
\hfill$\Box$

\begin{remark}
Not for every cubical set $X$ the sequence $X^N_n$ will form a semi-cubical set. Therefore, it cannot be said that Theorem \ref{comloc} states that homology with coefficients in a local system is isomorphic to the homology of some semi-cubical set.
\end{remark}

\subsection{Cohomology of cubical sets with coefficients in local systems}

Let $\mA$ be an AB4*-category. Contravariant systems of objects in the category 
$\mA^{op}$ are called covariant systems in $\mA$. Colimits in $\mA^{op}$ will be limits in $\mA$. Thus, a {\it covariant system} of Abelian groups on a cubical set $X$ is a diagram $F: \Box/X\to \Ab$.
Theorem \ref{main2} implies that for any cubical set $X$ the cohomology groups $H^n(X,F)$ with coefficients in the covariant system $F: \Box/X\to \Ab$ are isomorphic 
to $\Lim^n_{\Box/X}F$.
Covariant and contravariant systems are called {\it local} if their values on morphisms of the category $\Box/X$ are isomorphisms.

There is a statement proved by Quillen \cite[Chapter II, \S3, Proposition 4]{qui1967} characterizing weak equivalences in the category of simplicial sets.
 It follows from this assertion that for any weak equivalence of simplicial sets  $f: X\to Y$ and for a local system $F: \Delta/Y \to \Ab$, the cohomology group homomorphisms $H^n(Y, F) \to H^n(X, F\circ (\Delta/f))$ are isomorphisms for all $n\geq 0$.

Local systems on a simplicial set are defined by Gabriel and Zisman \cite[Appendix II, \S4.5]{gab1967} as contravariant systems $F: (\Delta/X)^{op}\to \mA$ taking values in an arbitrary abelian category $\mA$ with exact coproducts. Homologies are defined as satellites of the colimit functor.
If $\mA=\Ab^{op}$ is substituted, then the homology becomes the cohomology $H^n(X,F)$ of simplicial sets $X$ with coefficients in local systems $F$, which, according to the definition of Gabriel and Zisman, will be equal to $\Lim^n_{\Delta/X}F$.

The question arises whether the homomorphisms of the cohomology groups $H^n(Y, F)\to H^n(X, F\circ (\Box/f))$ with coefficients in the local system $F$ are isomorphisms in the case when $f : X\to Y$ is a weak equivalence of cubical sets.

Recall the definition of a test category.
Denote by $\Hot$ the classical homotopy category constructed from the category of topological spaces or simplicial sets as the category of fractions with respect to weak equivalences.
 
 It is known \cite{gro1983} that the category $\Hot$ is equivalent to the category of fractions for $\Cat$
 with respect to the class $W_{\infty}$ consisting of arrows $f: \mC\to \mC'$ of the category $\Cat$, for which the simplicial mapping of nerves $\nr(f): \nr\mC\to \nr \mC'$ is a weak equivalence (of simplicial sets).
 The equivalence of the categories $\Hot \to W^{-1}_{\infty} \Cat$ is constructed using the functor $Simpl: \Set^{\Delta^{op}}\to \Cat$
(see \cite[page 5]{mal2005})
assigning to each simplicial set the category of its simplices. For any $X\in \Set^{\Delta^{op}}$, this category is isomorphic 
to $\Delta/X$.
The category $\Hot$ is equivalent to the category of fractions of the category of simplicial sets with respect to morphisms $f: X\to Y$ such that $\Delta/f: \Delta/X \to \Delta/Y$ is a weak equivalence in $ \Cat$.

Problems related to test categories are briefly and clearly described in \cite{buc2017}, devoted to classes of test categories of cubical sets.
  
 Let $\mD$ be a small category. Consider the functor $\mD/(-): \mD\to \Cat$, assigning to each object $a\in \mD$ the category $\mD/a$, and to each morphism $\alpha: a\to b$ the functor $\mD/\alpha: \mD/a \to \mD/b$, taking every object $a'\to a$ of the category $\mD/a$ to the composition $a'\to a \xrightarrow{\alpha} b$.
The following diagram corresponds to the functor $\mD/(-)$
$$
 \xymatrix{
 \Set^{\mD^{op}} \ar[rr]^{i_{\mD}} && \Cat\\
 & \mD\ar[ru]_{\mD/(-)} \ar[lu]^{h^{\mD}}
 }
 $$ 
  where the functor $i_{\mD}$ is the left Kan extension of the functor $\mD/(-)$ along the Yoneda embedding $h^{\mD}$. It takes the values $i_{\mD}(X)= \mD/X$ on objects of the category $\Set^{\mD^{op}}$, and $i_{\mD}(f) = \mD/f $ - on morphisms.
  The functor $i_{\mD}$ has a right adjoint functor $i^*_{\mD}: \Cat\to \Set^{\mD^{op}}$.

On the category of presheaves $\Set^{\mD^{op}}$, consider the class $W_{\mD}$ consisting of morphisms of presheaves $f: X\to Y$ for which mappings of nerves $ \nr\mD/f: \nr\mD/X \to \nr\mD/Y$ are weak equivalences.
   It is clear that $W_{\mD}= i^{-1}_{\mD}(W_{\infty})$.
  
   Since the functor $i_{\mD}$ maps morphisms from $W_{\mD}$ to morphisms from $W_{\infty}$, it induces a functor between categories of fractions
    $\overline{i_{\mD}}: W^{-1}_{\mD}\Set^{\mD^{op}}\to W^{-1}_{\infty} \Cat \simeq \Hot$. 
    
The following definition of test categories is given in \cite{buc2017}:

\begin{enumerate}
\item A category $\mD$ is called a weak test category if the functor $\overline{i_{\mD}}$ is an equivalence of categories.

\item A category $\mD$ is called local test category if for each $a\in \Ob\mD$ the comma category $\mD/a$ is weak test category.

\item A category $\mD$ is called test if it is weak test and locally test.

\item A category $\mD$ is strict test if it is test category and if the functor $i_{\mD}: \Set^{\mD^{op}} \to \Hot$ commutes with finite products.
\end{enumerate}

If $\mD$ is a test category, then the category $\Set^{\mD^{op}}$ has a closed model structure in which cofibrations are monomorphisms and weak equivalences are all morphisms belonging to the class $W_{\mD}$. It is clear that in this case there is an equivalence
homotopy categories $W^{-1}\Set^{\mD^{op}}\simeq \Hot$.
This model structure is called the standard \cite{jar2006}.

In \cite{cis2002} it was proved that the category of cubical sets is a test category.
This implies that the category of cubical sets has a standard model structure. Weak equivalences of cubical sets with respect to this model structure are morphisms $f: X\to Y$ for which the corresponding simplicial mapping
nerves of the categories ${\nr}(\Box/f): {\nr}(\Box/X)\to {\nr}(\Box/Y)$ is a weak equivalence of simplicial sets.

\begin{corollary}\label{lociso}
Let $f: X\to Y$ be a weak equivalence of cubical sets with respect to the standard model structure.
Then for any local system $L: \Box/Y\to \Ab$ and $n\geq 0$ the homomorphism $H^n(Y, L)\to H^n(X,L\circ(\Box/f))$ is an isomorphism.
\end{corollary}
{\sc Proof.} Let $\partial: \Delta/{\nr}\mC \to \mC$ denote the functor that assigns to each simplex $c_0\to \cdots \to c_n$ the object $c_n\in \mC$.
As in the example \ref{hsimp}, $H^{simp}_n$ denote the simplicial homology of the category nerve, for all $n\geq 0$.
It is well known that the functor $\partial$ is aspherical, and hence $H^{simp}_n(\partial/c)\cong H^{simp}_n(\pt)$,
whence for each functor $G: \mC\to \Ab$ there is a natural isomorphism $\lim^n_{\mC}G \stackrel{\cong}\to \lim^n_{\mC}G\partial$.
Consider the commutative diagram

$$
\xymatrix{
\Delta/ {\nr}(\Box/X)\ar[d]_{\partial} \ar[rr]^{\Delta/{\nr}(\Box/f)} &&
\Delta/ {\nr}(\Box/Y)\ar[d]_{\partial} \ar@{-->}[rd]^{L\partial} \\
\Box/X \ar[rr]_{\Box/f} && \Box/Y \ar[r]_L & \Ab
}
$$
It follows from  \cite[Chapter II, \S3, Proposition 4]{qui1967} that for any weak equivalence of simplicial sets, the corresponding homomorphism of cohomology groups with coefficients in a local system is an isomorphism.
Since, by assumption, the simplicial mapping ${\nr}(\Box/f): {\nr}(\Box/X)\to {\nr}(\Box/Y)$ is a weak equivalence, we have isomorphism
$\Lim^n L\partial \stackrel{\cong}\to
\Lim^n L\partial(\Delta/{\nr}(\Box/f))$.
Commutative diagram of natural homomorphisms
$$
\xymatrix{
\Lim^n L (\Box/f)\partial & \Lim^n L\partial(\Delta/{\nr}(\Box/f)) \ar[l]_{=}&
\ar[l]_{\qquad\cong} \Lim^n L\partial\\
\Lim^n_{\Box/X}L\circ(\Box/f)\ar[u]_{\cong} && \Lim^n_{\Box/Y}L
\ar[ll]\ar[u]_{\cong}
}
$$
shows that the homomorphism in the bottom row is an isomorphism.
By Theorem \ref{main2}, this implies that the canonical homomorphism $H^n(Y,L)\to H^n(X, L\circ(\Box/f))$ is an isomorphism
 for all $n\geq 0$.
 \hfill $\Box$

\begin{example}
Morphisms between standard cubes do not preserve cohomology groups.
For example, for a unique morphism from a standard one-dimensional cube to a cubical point $s=\Box(\sigma^1_1): \Box^1\to \Box^0$ the inverse fiber $\overleftarrow{s}(y)$ of the cube $\tilde{y}= \Box(\sigma^1_1): \Box^1\to \Box^0$ is isomorphic to the cubical
 set $\Box^1 \times\Box^1$ having nonzero homology groups 
 $H_2(\Box^ 1\times\Box^1, \ZZ)\cong H_1(\Box^1 \times\Box^1, \ZZ)=\ZZ$ \cite[Example 1]{X2019}. This means that there are covariant systems 
 $F: \Box/\Box^0 \to \Ab$ and a number $n>0$ such that the group homomorphism $H^n(\Box^1, s^* F) \to H^ n(\Box^0, F)$ is not an isomorphism.

The cubical set morphism $s=\Box(\sigma^1_1): \Box^1\to \Box^0$ 
is a weak equivalence.
 For any local system $F$ on the cubical point $\Box^0$, $H^n(\Box^0, F)$ and $H^n(\Box^1, F)$ will be equal to $0$ for $ n>0$, and the homomorphism $H^0(\Box^0,F)\to H^0(\Box^1,s^* F)$ is an isomorphism.
\end{example}

\begin{remark}
For an arbitrary small category $\mD$, we can consider the category $\Set^{\mD^{op}}$ as a category with weak equivalences, in the sense of \cite{cha2001}. We can declare any morphism $f: X\to Y$ to be a weak equivalence if the induced simplicial mapping ${\nr}(\mD/X)\to {\nr}(\mD/Y)$ is a weak equivalence of simplicial sets.
If we take as a definition $H^n(X,L):= \Lim^n_{\mD/X} L$ the cohomology group of
 the $\mD$-set $X\in \Set^{\ mD ^{op}}$ with coefficients in the local system 
 $L: \mD/X\to \Ab$, then the corollary \ref{lociso} and its proof remain true
  if the category of cubes is replaced by $\mD$.
\end{remark}

\subsection{Spectral sequence for a morphism of cubical sets}
 
 Theorem \ref{main2} allows us to apply Proposition \ref{spmor}, Corollary \ref{lociso}, and Remark \ref{spmorinv} to construct a spectral sequence for a morphism of cubical sets endowed with local systems. This spectral sequence links the cohomology groups of cubical sets and the cohomology of inverse fibers.
Weak equivalences of cubical sets belonging to the standard model structure are considered.

\begin{corollary}\label{spseqserr}
Let $f: X\to Y$ be a cubical set morphism whose inverse fiber diagram morphisms are weak equivalences. Then for any local system of Abelian groups $G: \Box/Y \to \Ab$ there exists a spectral sequence of the first quarter
$$
E^{pq}_2 = H^p\left( Y, 
\{H^q( \overleftarrow{f}(\sigma), f^*_{\sigma}G)\}^{-1}_{\sigma\in \Box/Y}\right)
\Rightarrow H^{p+q}(X, G),
$$
where $f^*_{\sigma}G$ equals the composion of 
$\Box/\overleftarrow{f}(\sigma)\xrightarrow{\Box/f_{\sigma}} 
\Box/X \xrightarrow{G} \Ab$.
\end{corollary}

\section{Homology of semi-cubical sets as cubical homology}

Recall that the category $\Box_+$ has a set of objects $\Ob\Box_+= \Ob\Box$.
Its morphisms are all monomorphisms of the cube category $\Box$.
A semi-cubical set is a functor $\Box^{op}_+\to \Set$.
In this section, we prove that the homology of a semi-cubical set with coefficients in a contravariant system is isomorphic to the homology of the universal cubical set containing it.

Note that the non-degenerate cubes of the universal cubical set coincide with the cubes
 of the given semi-cubical set, and the values of the contravariant system on the non-degenerate cubes of the universal cubical set are equal to the values 
 of the contravariant system on the corresponding cubes of the given semi-cubical set.

\subsection{Universal cubical set for a given semi-cubical set}

Let $J: \Box_+\to \Box$ be an embedding functor.
The next lemma shows that this functor is a virtual discrete prefibration.

\begin{lemma}\label{confplus}
For any $\II^n\in \Ob\Box$ each connected component of the category $\II^n/J$ has an initial object.
\end{lemma}

\noindent{\sc Proof.} 
Objects of the category $\II^n/J$ can be considered as
morphisms $\II^n\to \II^k$ of the category $\Box$. The morphism between $\II^n\stackrel{\alpha}\to \II^p$ and $\II^n \stackrel{\beta}\to \II^q$ in this category is given by the commutative triangle
\begin{equation}\label{morphcomma}
\xymatrix{
& \II^n \ar[ld]_{\alpha} \ar[rd]^{\beta} \\
\II^p \ar@{>-->}[rr]^{\mu} && \II^q
}
\end{equation}
where $\mu$ is a monomorphism of the category $\Box$ of cubes.
Denote this morphism of the category $\II^n/J$ by
$\alpha \stackrel{\mu}\to \beta$.
Each morphism $\II^n \stackrel{\beta}\to \II^q$ of the category of cubes admits a unique decomposition $\beta= \mu\circ \sigma$ into the composition of an epimorphism $\sigma$ and a monomorphism $\mu$.
Denote the epimorphism $\sigma$ by $i(\beta)$ and the monomorphism $\mu$ by $i_{\beta}$. We obtain a morphism $i(\beta)\stackrel{i_{\beta}}\to \beta$, and the set of morphisms $i(\beta)\to \beta$ in the category $\II^n/J$ consists of a single element equal to $i_{\beta}$.

Expanding triangle morphisms (\ref{morphcomma}) into compositions of epimorphisms and monomorphisms in the category of cubes, we obtain a commutative diagram

$$
\xymatrix{
\II^k \ar@{>->}[rd]_{i_{\alpha}} && \II^n \ar@{->>}[ll]_{i(\alpha)} \ar[ld]_{\alpha}
\ar[rd]^{\beta} \ar@{->>}[rr]^{i(\beta)}
&& \II^m \ar@{>->}[ld]^{i_{\beta}} \\
& \II^p \ar@{>->}[rr]^{\mu} && \II^q
}
$$
Here, epimorphisms are represented by the arrow $\twoheadrightarrow$, and monomorphisms by the arrow $\rightarrowtail$.
Since the decomposition of the morphism $\beta$ into the composition of an epimorphism and a monomorphism is unique, this diagram leads to the equalities:

$$
i(\alpha)= i(\beta), \quad \II^k=\II^m, \quad \mu\circ i_{\alpha}=i_{\beta}.
$$
We arrive at the following commutative diagram in the category
 $\II^n/J$:

$$
\xymatrix{
i(\alpha)\ar@{=}[r] \ar[d]_{i_{\alpha}} & i(\beta)\ar[d]^{i_{\beta}}\\
\alpha \ar[r]_{\mu} & \beta
}
$$
This implies that all objects $\alpha$ of the category $\II^n/J$ belonging to the same connected component have the same object $i(\alpha)$, which is an epimorphism in the category of cubes.

This object is the initial object of the connected component containing the object $\alpha$ and the morphism $i_{\alpha}: i(\alpha)\to \alpha$ is the unique morphism from this object to the object $\alpha$.
\hfill$\Box$

\medskip

Let us proceed to the construction of a universal cubical set for a semi-cubical set.
As in the proof of Lemma \ref{confplus},
for any morphism $\alpha: \II^m\to \II^n$ we denote by $\II^m\stackrel{i(\alpha)}\to \II^p$ the epimorphism, and by $\II^p\stackrel{i_{\alpha}}\to \II^n$ monomorphism of the cube category such that $\alpha= i_{\alpha}\circ i(\alpha)$. For any family of sets $(S_i)_{i\in I}$ its disjoint union $\coprod\limits_{i\in I}S_i$ consists of pairs $(i,s)$ such that $s\in S_i$.

\begin{proposition}\label{cubeforsemi}
Let $X: \Box_+^{op}\to Set$ be an arbitrary semi-cubical set.
Then its left Kan extension $Lan^{J^{op}}X: \Box^{op}\to \Set$ is a functor assigning to each cube $\II^n$ a set
$$
(Lan^{J^{op}}X)_n=
\coprod\limits_{\II^n\stackrel{\gamma}\twoheadrightarrow\II^k}X_k,
$$
and to each morphism $\alpha: \II^m\to \II^n$ a mapping of sets
$$
Lan^{J^{op}}\alpha=\alpha_*:
\coprod\limits_{\II^n\stackrel{\gamma}\twoheadrightarrow\II^k}X_k
\to
\coprod\limits_{\II^m\stackrel{\beta}\twoheadrightarrow\II^p}X_p,
$$
defined as
$\alpha_*(\gamma,x)= (i(\gamma\circ\alpha), X(i_{\gamma\circ\alpha})(x))$.
\end{proposition}
{\sc Proof.}
Consider the functor $J: \Box_+\to \Box$. By Lemma \ref{confplus}, for every object 
$\II^n\in \Box$ each connected component of the category $\II^n/J$ has an initial object. This allows us to use the construction from Proposition \ref{lemma36} to construct the functor $Lan^{J^{op}}X$. This functor on objects $\II^n$ takes the values $\coprod\limits_{\gamma\in init(\II^n/J)} X Q^{op}_{\II^n}$, where $init( \II^n/J)$ is the set of initial objects in the 
connection components of the category $\II^n/J$, and $Q_{\II^n}: \II^n/J\to \Box_+$
 is a functor that assigns to each object $\II^n\to J(\II ^p)$ of the category $\II^n/J$ the object $\II^p\in \Box_+$, and the morphism

$$
\xymatrix{
& \II^n \ar[ld] \ar[rd] \\
J(\II^p) \ar[rr]_{J(\delta)} && J(\II^q)
}
$$
maps the morphism
 $\II^p\stackrel{\delta}\to \II^q$ in the category $\Box_+$.
This implies that
$Lan^{J^{op}}X (\II^n)=
 \coprod\limits_{\II^n\stackrel{\gamma}\twoheadrightarrow \II^k}X_k$.
Applying Proposition \ref{lemma36} clause leads to the definition$\alpha_*$ 
as a mapping of sets, making the diagram commutative
$$
\xymatrix{
\coprod\limits_{\II^n\stackrel{\gamma}\twoheadrightarrow \II^k}X_k
\ar@{-->}[rr]^{\alpha_*} &&
\coprod\limits_{\II^m\stackrel{\beta}\twoheadrightarrow \II^p}X_p\\
X_k \ar[rr]_{X(i_{\gamma\circ\alpha})}\ar[u]^(.35){in_{\gamma}}
&& X_p\ar[u]_(.35){in_{i(\gamma\circ\alpha)}}
}
$$
where $in_{\gamma}$ are mappings defined as $in_{\gamma}(x)=(\gamma,x)$. This diagram leads to the desired formula for the mapping $\alpha_*$.
\hfill$\Box$

\begin{corollary}
Cube $(\beta,y)\in (Lan^{J^{op}}X)_m$ consisting of 
epimorphism $\beta: \II^m\to \II^p$ and element $y\in X_p$
is degenerate if and only if $m\neq p$.
\end{corollary}
{\sc Proof.} 
The cube $(\beta,y)$ is degenerate if and only if there exist an epimorphism $\alpha: \II^m\to \II^n$, where $n=m-1$, and an element $(\gamma,x )\in
\coprod\limits_{\II^n\stackrel{\gamma}\twoheadrightarrow \II^k}X_k$ such that $\alpha_*(\gamma, x)= (\beta,y)$.
The epimorphism $\alpha$ will be equal to the degeneration map $\sigma^m_i$ for some $i$.
The last equality is equivalent to the condition $(i(\gamma\circ\alpha), X(i_{\gamma\circ\alpha})(x))= (\beta,y)$.
Since $i_{\gamma\alpha}\circ i(\gamma\alpha)$ is a canonical decomposition of $\gamma\alpha$ into a composition of an epimorphism and a monomorphism, and $\gamma\alpha$ is an epimorphism, then $i(\gamma\alpha)= \gamma\alpha$ and $i_{\gamma\alpha}=1_{\II^k}$. We conclude that $(\beta,y)$ is degenerate if and only if $\beta=\gamma\alpha$ and $y=x$, and $k=p$, $n=m-1$. Thus, if $m>p$, then
$(\beta,y)= \alpha_*(\gamma,y)$, where $\alpha=\sigma^m_i: \II^m\to \II^{m-1}$, for some $i\in\{1, \ldots, m\}$.
\hfill$\Box$

\medskip

For every $n\geq 0$, the set of $n$-dimensional non-degenerate cubes of the cubical set $Lan^{J^{op}}X$ is equal to $\{1_{\II^n}\}\times X_n$.
Applying Proposition \ref{cubeforsemi} to $\alpha= \delta^{n,\varepsilon}_i$, we conclude that the boundary operator $(\delta^{n,\varepsilon}_i)_*$ transports a non-degenerate cube $(1_{\II^n},x)$ into a non-degenerate cube $(1_{\II^{n-1}}, 
X(\delta^{n, \varepsilon}_i))$.
We get

\begin{corollary}
The adjunction unit $\eta_X: X\to (Lan^{J^{op}}X)\circ J^{op}$
defines an isomorphism between the semi-cubical set $X$ and the semi-cubical set consisting of 
non-degenerate cubes and boundary operators of the cubical set $Lan^{J^{op}}X$.
\end{corollary}

\subsection{Comparison of homology for cubical and semi-cubical sets}

Let $X: \Box^{op}_+\to Set$ be a semi-cubic set, $Lan^{J^{op}}X: \Box^{op}\to Set$ be the corresponding cubic set.
Before comparing cubical and semi-cubical homology, we first study the structure of the category $\Box/Lan^{J^{op}}X$, using Proposition \ref{cubeforsemi}.

For this purpose we denote $Y= Lan^{J^{op}}X$.
Since $Y_n= \coprod\limits_{\gamma:\II^n\twoheadrightarrow \II^k}X_k$, then the objects of the category $\Box/Lan^{J^{op}}X$ consist of pairs $(\gamma,x)$, where $\gamma: \II^n\to \II^k$ is an epimorphism in the category of cubes, and $x\in X_k$ is an arbitrary element of the semi-cubical set $X$.
They can be considered as pairs of morphisms in the category $\Set^{\Box^{op}}$:
$$
  \Box^n\xrightarrow{\Box{(\gamma)}} \Box^k\stackrel{\widetilde{x}}\to X,
$$
where $\gamma$ is an epimorphism in the category $\Box$ and $x$ is an element in 
$X_k$.
If $k=n$ and $\gamma=1_{\II^n}$, then this pair is a non-degenerate cube 
of the cubical set $Lan^{J^{op}}X$.
Morphisms of the category $\Box/Lan^{J^{op}}X$ are defined using commutative triangles

$$
\xymatrix{
& Y  \\
\Box^m \ar[ru]^{\widetilde{(\beta,y)}} \ar[rr]_{\Box(\alpha)}
&& \Box^n \ar[lu]_{\widetilde{(\gamma,x)}}
}
$$
meaning that the mapping $\alpha_*: Y_n\to Y_m$, defined in Proposition \ref{cubeforsemi}, 
satisfies the condition $\alpha_*(\gamma,x)= (\beta,y )$.
This condition is equivalent to two equalities $i(\gamma\circ \alpha)=\beta$ and $X(i_{\gamma\circ\alpha})(x)=y$.

This implies that the morphism $\widetilde{(\beta,y)}\to \widetilde{(\gamma,x)}$ can be considered as a commutative diagram

\begin{equation}\label{morphlan}
\xymatrix{
\Box^m\ar[d]_{\Box(\alpha)} \ar[r]^{\Box(\beta)}
& \Box^p \ar[d]_{\Box(\mu)} \ar[r]^{\widetilde{y}}& X\\
\Box^n \ar[r]_{\Box(\gamma)} & \Box^k \ar[ru]_{\widetilde{x}}
}
\end{equation}
in which $\mu$ is a monomorphism and $\beta$ and $\gamma$ are epimorphisms in the category of cubes. The commutativity of this diagram is equivalent to the three equalities $i(\gamma\circ\alpha)=\beta$, $i_{\gamma\circ\alpha}=\mu$, $\widetilde{x}\Box(\mu)= \widetilde{y}$, which are equivalent to $\alpha_*(\gamma,x)=(\beta,y)$.

We have proved the following assertion.

\begin{proposition}
Each object of the category $\Box/Lan^{J^{op}}X$ can be given by a pair of morphisms 
$(\II^n \stackrel{\gamma}\to \II^k,\Box^k \stackrel{\widetilde{x}}\to X)$ where 
$\gamma$ is an epimorphism in the cube category.
Morphisms $(\II^m\stackrel{\beta}\to \II^p,\Box^p
\stackrel{\widetilde{y}}\to X) \to(\II^n\stackrel{\gamma}\to \II^k,\Box^k \stackrel{\widetilde{x}}\to X) $, are pairs $(\alpha,\mu)$ consisting of the morphism $\II^m\stackrel{\alpha}\to \II^n$ and the monomorphism
$\II^p\stackrel{\mu}\to \II^k$ in the cube category, making the diagram (\ref{morphlan}) commutative.
\end{proposition}

To compare homology, we need one more statement:

\begin{proposition}\label{semcubadj}
Functor $J_*: \Box_+/X \to \Box/Lan^{J^{op}}X$ assigning each $\widetilde{x}: \Box^n_+\to X$ cube $\widetilde {(1_{\II^n},x)}$ has a left adjoint $S: \Box/ Lan^{J^{op}}X \to \Box_+/X$ defined on objects as $S( \gamma, \widetilde{x})= \widetilde{x}$, and on morphisms $S(\alpha,\mu)= \mu$.
\end{proposition}
{\sc Proof.} If we stick to the definition of objects and morphisms of the category $\Box/Lan^{J^{op}}X$ shown by the diagram (\ref{morphlan}), this will greatly simplify the proof.
The assertion will follow from the universality
arrows $\widetilde{(\gamma, x)}\xrightarrow{(\gamma, \Box(1_{\II^k}))} J_*(\tilde{x})= (1_{\II^k} , \tilde{x})$,
for every $\widetilde{(\gamma, x)}\in \Box/ Lan^{J^{op}}X$. The adjointness
 of functors follows from \cite[Theorem IV.1.2, Page 83]{mac1972}.
\hfill$\Box$

\medskip

The paper \cite{X2008} gives a definition for the homology groups of a semi-cubical set.
Similarly, one can define the homology of a semi-cubical set with coefficients in the contravariant system of objects $F$ as the homology of the complex consisting of objects $C_n(F)= F(\II^n)$ and differentials
$d_n= \sum\limits^n_{i=1}(-1)^i (F(\delta^{n,0}_i)- F(\delta^{n,1}_i))$.
For any semi-cubical object $F: \Box^{op}_+\to \mA$ the tensor product of $F$ and the projective resolution of the functor $\Delta_{\Box_+}{\mathbb Z}$ used in the proof \cite[Proposition 4.2]{X2008} is equal to this complex $C_*(F)$. This leads to a generalization of the assertion \cite[Theorem 4.3]{X2008} and gives an isomorphism between the values $\coLim_n^{(\Box_+/X)^{op}}F$ and the homology objects of the semi-cubical set $X$ with coefficients in $F$, which we denote by $H^+_n(X,F)$.
The next assertion shows that they are isomorphic to the homology of the universal cubical set containing $X$ with coefficients in $F\circ S^{op}$.

\begin{theorem}\label{semicubecube}
Let $X$ be a semi-cubical set, and $\mA$ an AB4-category. Consider the functor $S: \Box/ Lan^{J^{op}}X\to \Box_+/X$ taking each cube $(\gamma,\widetilde{x})$ to the cube $\widetilde{x} $, and the morphism $(\alpha,\mu)$ into the morphism $\mu$.
Then for every functor $F: (\Box_+/X)^{op}\to \mA$
there are natural isomorphisms
$H_n(Lan^{J^{op}}X, FS^{op})\to H^+_n(X,F)$,
for all $n\geq 0$.
\end{theorem}
{\sc Proof.}
By Theorem \ref{main2},  for all $n\geq 0$, the objects
$H_n(Lan^{J^{op}}X, FS^{op})$ are isomorphic to 
$\coLim_n^{(\Box/Lan^{J^{op} }X)^{op}} FS^{op}$.
Since $S$ has a right adjoint, the categories $S/\widetilde{x}$ have point homology, 
and hence, according to \cite[Theorem 2.3]{obe1968}, the canonical morphisms
$$
\coLim_n^{(\Box/Lan^{J^{op}}X)^{op}} FS^{op}\to \coLim_n^{(\Box_+/X)^{op}}F
$$
are isomorphisms. This implies the existence of the required isomorphism.
\hfill$\Box$

\section{Conclusion}

In this paper, we study the homology of cubical sets with coefficients in contravariant systems whose morphisms must not be isomorphisms. The following facts have been established:
\begin{itemize}
\item These homology are invariant under morphism between cubical sets when passing
 to the direct image of the system of coefficients.
\item There is a criterion for the invariance of these homologies when passing to the inverse image.
\item
These homology generalize the singular cubical homology with local coefficients and the homology of semi-cubical sets
with coefficients in contravariant systems.
\item
There is a spectral sequence for colimit homologies of cubical sets with coefficients in 
contravariant systems.
\item
There is a spectral sequence for morphism domain cohomology between cubical sets with local systems of Abelian groups.
\item The homology of small category with coefficients in a diagram can be calculated as cubical homology.
\item
The Baues-Wirsching cohomologies with coefficients in natural systems are isomorphic to cubical cohomologies with coefficients in covariant systems.
\end{itemize}

We hope that these results will find applications in solving problems related to the methods of calculating homological groups for cubical sets, in studying homology for mathematical models of computing systems and processes, in studying the homology of topological spaces with local coefficients.

\end{document}